\input ./style/arxiv-general.cfg
\documentclass[MSNbibl,number,citesort,seceqn,dvips]{arxbj}
\makeatletter
   \@ifpackageloaded{graphicx}{}{\usepackage{graphicx}}
\makeatother
\usepackage{dcolumn}

\aid{0}
\volume{21}
\issue{4}
\pubyear{2015}
\firstpage{2217}
\lastpage{2241}
\doi{10.3150/14-BEJ642} 
\docsubty{FLA}

\makeatletter

\newcommand{\lleft}{\left}
\newcommand{\rright}{\right}
\newtheorem{lea}{Lemma}[section]
\newtheorem{cor}{Corollary}[section]
\newtheorem{prn}{Proposition}[section]
\newremark{rem}{Remark}[section]
\newproclaim{Assumptions}{Assumptions}

\newcommand{\UE}{\operatorname{UE}}
\renewcommand{\Re}{\operatorname{Re}}
\newcommand{\Var}{\operatorname{Var}}
\newcommand{\rank}{\operatorname{rank}}
\newcommand{\diag}{\operatorname{diag}}
\newcommand{\E}{\mathrm{E}}

\newcommand{\eqref}[1]{(\ref{#1})}
\newcommand{\xrightarrow}{\hbox to 30pt{\rightarrowfill}}

\newcommand{\ADR}{\operatorname{ADR}}
\newcommand{\trace}{\operatorname{trace}}
\newcommand{\RMSE}{\operatorname{RMSE}}
\newcommand{\risk}{\operatorname{risk}}

\renewcommand{\epsilon}{\varepsilon}
\renewcommand{\leqslant}{\leq}
\renewcommand{\geqslant}{\geq}
\newcommand{\SSR}{\operatorname{SSR}}
\newcommand{\Ch}{\operatorname{Ch}}

\newcolumntype{d}[1]{D{.}{.}{#1}}

\makeatother

\begin{document}
\begin{frontmatter}

\title{Optimal method in multiple regression with structural changes}
\runtitle{Stein rules in linear models with change-points}

\begin{aug}
\author[A]{\inits{F.}\fnms{Fuqi}~\snm{Chen}\thanksref{e1}\ead[label=e1,mark]{chen111n@uwindsor.ca}} \and
\author[A]{\inits{S.}\fnms{S\'ev\'erien}~\snm{Nkurunziza}\corref{}\thanksref{e2}\ead[label=e2,mark]{severien@uwindsor.ca}}
\address[A]{Mathematics and Statistics Department, University of
Windsor, 401 Sunset Avenue, Windsor,
Ontario N9B 3P4, Canada. \printead{e1,e2}}
\end{aug}

\received{\smonth{10} \syear{2013}}
\revised{\smonth{5} \syear{2014}}

%
\begin{abstract}
In this paper, we consider an estimation problem
of the regression coefficients in multiple regression models
with several unknown change-points. Under some realistic
assumptions, we propose a class of estimators which includes as a
special cases shrinkage estimators (SEs) as well as the unrestricted
estimator (UE) and the restricted estimator (RE). We also derive a
more general condition for the SEs to dominate the UE. To this end,
we generalize some identities for the evaluation of the bias and
risk functions of shrinkage-type estimators. As illustrative
example, our method is applied to the ``gross domestic product''
data set of 10 countries whose USA, Canada, UK, France and Germany.
The simulation results corroborate our theoretical findings.
\end{abstract}

%
\end{frontmatter}

\section{Introduction}
In this paper, we study the multivariate regression models with
multiple change-points occurring at unknown times. The target
parameters are the regression coefficients while the unknown change
points are treated as nuisance parameters. More specifically, we are
interested in scenario where imprecise prior information about the
regression coefficients is available, that is, the target parameters may
satisfy some restrictions.

The importance of change-points' model in literature is a primary
source of our motivation. Indeed, the regression model with
change-points has been applied in many fields. For example, this
model was used in Broemeling and Tsurumi \cite{bt1987} for the US demand
for money, as well as in Lombard \cite{l1986} for the
effect of sudden
changes in wind direction of the flight of a projectile. It was also
analyse the DNA sequences (see, e.g., Braun and Muller \cite
{bm1998} and Fu
and Curnow \cite{fc1990a,fc1990b}). To give some
recent references, we quote
Bai and Perron \cite{bp2003}, Zeileis \emph{et~al.} \cite{zei},
Perron and Qu \cite{pq2006} among others.

More specifically, the method in Perron and Qu \cite{pq2006} is
based on a
global least squares procedure. Generally, when the restriction
holds, the restricted estimator (RE) dominates the unrestricted
estimator (UE). However, it is well known that the RE may performs
poorly when the restrictions is seriously violated.

Over the years, shrinkage estimation has become a useful tool in
deriving the method which combines in optimal way both imprecise
prior knowledge from a hypothesized restriction and the sample
information. For more details about such a technique, we refer to
James and Stein \cite{js1961}, Baranchick \cite
{baranchick1964}, Judge and Bock \cite{judg1978},
and the references therein. Also, to give some recent contributions
about shrinkage methods, we quote Saleh \cite{sal},
Nkurunziza and Ahmed \cite{nku}, Nkurunziza \cite{n2011} and {
Tan \cite{tan}}, among others.

To the best of our knowledge, in context of multiple regression
model with unknown changes-points, shrinkage method has received, so
far, less attention. Thus, we hope to fill this gap by developing a
class of shrinkage-type estimators which includes as special cases
the UE, RE, James--Stein type and positive shrinkage estimators as
well as pre-test estimators. We also prove that the proposed
shrinkage estimators (SEs) dominate in mean square error sense the
UE. The technique in this paper extends, in two ways the method
given in literature.

First, the asymptotic dependance structure between the shrinking
factor (i.e., the difference between the UE and the RE) and the RE is
more general than that given in the quoted papers. In particular,
the asymptotic variance of RE and the asymptotic variance of
$(\UE-\Re)$ are not positive definite matrices as in
the problem studied in Judge and Mittelhammer \cite{judg2004}. { This is
justified by the fact that, since the hypothesized restriction is
linear, these quantities are asymptotically equivalent to the
nonsurjective linear (equivalent here to noninjective linear)
transformations of the UE for which the asymptotic variance is
positive definite matrix. In this case, it is impossible for the
asymptotic variance of RE or that of $(\UE-\Re)$ to
be positive definite matrix. To make the justification more precise,
let $\mathbf{A}$ be a nonrandom $n\times m$-matrix with the rank
$n_{0}<n$, let $B$ be a nonrandom $n$-column vector, and let $F$ be
$n$-column random vector whose variance is a positive definite
matrix $\boldsymbol{\Psi}$. Further, let $G=A F+B$, that is a nonsurjective
linear transformation of the random vector $F$. Then,}
{$\Var(G)=\mathbf{A}\boldsymbol{\Psi} \mathbf{A}'$ which cannot be a
positive definite matrix since $\rank(\mathbf{A}\boldsymbol{\Psi}
\mathbf{A}')=n_{0}<n$.}

Second, we derive a more general condition for the SEs to dominate
the UE. To this end, we generalize Theorem~1 and Theorem~2 of
Judge and Bock \cite{judg1978} which are useful in computing
the bias and the
risk functions of shrinkage-type estimators. As far as the
underlying asymptotic results are concerned, another difference,
with the work in Judge and Mittelhammer \cite{judg2004}, consists in
the fact
that we derived the joint asymptotic normality under weaker
conditions than that in the quoted paper. {Indeed, in Judge and Mittelhammer \cite{judg2004}, the covariance--variance of the error terms
is a
scalar matrix (see the first paragraph of Section~2 in Judge and Mittelhammer \cite{judg2004}) and thus, the errors term are both homoscedastic
and uncorrelated.} {In addition, in the quoted paper, the regressors
are} {assumed nonrandom.} {In this paper, the errors term do not
need to be homoscedastic and/or uncorrelated, and they may also be
nonstationary stochastic processes. Further, the regressors may be
random and in addition, they may be correlated with the error terms.
In summary, the proposed method is applicable to the statistical
model with familiar regularity conditions as assumed in Judge and
Mittelhammer \cite{judg2004}, see the last sentence of
Section~2.4, as well
as in unfamiliar regularity conditions for which the dependance
structure of the errors and regressors terms is as weak as that of
mixingale array. The model considered here takes also an account for
the possibility of the change-points phenomenon and, because of
this, the derivation of the joint  asymptotic normality
between the UE and RE is mathematically challenging.} Moreover, the
established results extend that given for example in Perron and Qu~\cite{pq2006}.

{In concluding this introduction, note that due to the conditions
discussed above which} { are weaker than that in the literature, the
construction of shrinkage-type estimators cannot be obtained by
applying the results given in the quoted papers. Further, the
derivation of the asymptotic distributional risk (ADR) of shrinkage
estimators (SEs) is challenging and the instrumental identities in
Judge and Bock \cite{judg1978}, Theorems 1 and 2, are not
useful. This
motivated us to generalize these identities. This constitutes one of
the aspects of the main results which are significant in reflecting
the difference with the quoted works. The second aspect, of the main
results which is significant in reflecting the difference with the
quoted works, can be viewed from the fact that the established ADR
has some extra terms and the risk dominance condition of SEs looks
quite complicated.}

The rest of this paper is organized as follows.
Section~\ref{sec:statmodel} describes the statistical model and
outlines the proposed estimation strategies.
Section~\ref{sec:jointnorm} gives the joint asymptotic normality of
the unrestricted and restricted estimators. In
Section~\ref{sec:shrink}, we introduce a class of shrinkage-type of
estimators for the coefficients and derive its asymptotic
distribution risks. Section~\ref{sec:numerical} presents some
simulation studies and an illustrative analysis of a real data set.
Section~\ref{sec:conclu} gives some concluding remarks and, for the
convenience of the reader, technical proofs are given in the
\hyperref[app]{Appendix}.

\section{Statistical model and assumptions}\label{sec:statmodel}

In this section, we present the statistical model as well as the
main regularity conditions. {As mentioned above, in this paper, we
focus on the model with change-points. Nevertheless, the proposed
method is useful in linear model without change-points. In this last
case, the derivation of the joint asymptotic normality between the
RE and UE is not as mathematically involved as in case of the model
with change-points.}

{
\subsection{The linear model without change-points}} {We consider
the
multiple linear regression model with $T$ observations for which
the response is a $T$-column vector $Y=(y_1,\ldots,y_T)'$, the
regressors is a $T\times q_{0}$-matrix $\bar{Z}$, the regression
coefficients is a $q_{0}$-column vector $\delta$, and the errors
term is a $T$-column vector $u$. In particular, we have
let}
%
\begin{equation}
\label{modelprel} Y=\bar{Z}\delta+u.
\end{equation}
{Further, we consider the scenario where a prior knowledge about
$\delta$ exists with some uncertainty. More specifically, we
consider the case where $\delta$ is suspected to satisfy the
following} { restriction}
%
\begin{equation}
\label{resprel} R\delta=r,
\end{equation}
{where $R$ is a known $k\times q_{0}$-matrix with rank $k\leqslant
q_{0}$, and $r$ is a known $k$-column vector.} { Under some
regularities conditions on the error terms and the regressors, the
shrinkage estimator for the parameter $\delta$ is available in
literature. To give some references, we quote Saleh \cite{sal},
Hossain \emph{et al.} \cite{ahm2009} among others. The shrinkage estimators
given in the quoted papers are members of the class of shrinkage
estimators which is established in this paper. Further, the
established condition for the risk dominance of shrinkage estimators
is more general than that given for example, in Saleh \cite{sal},
Hossain \emph{et al.} \cite{ahm2009}.}

The { proposed methodology} { is applicable to the model in
\eqref{modelprel} and \eqref{resprel} } { provided that the
conditions on the error and regressors terms are such that, as $T$
tends to infinity,}
\begin{enumerate}[2.]
\item[1.] the matrices $T^{-1} \bar{Z}^{0\prime}\bar{Z}^0$ and
$T^{-1} (\bar{Z}^{0\prime}uu^{\prime}\bar{Z}^0)$
converge in probability to nonrandom { $q_{0}\times q_{0}$-positive
and definite matrices;} {
\item[2.]
$T^{-1/2}\bar{Z}^{0\prime}u$ converges in distribution to a Gaussian
random vector whose variance--covariance is the limit in probability
of $T^{-1} \bar{Z}^{0\prime}\bar{Z}^0$.}
\end{enumerate}
%
{ These two points are generally satisfied in classical regression
models where the error terms are homoscedastic and independent, with
linearly independent regressors. In the sequel, we consider a very
general model with change-points and heteroscedastic as well as
possibly correlated errors term. The assumptions of the model are
discussed in the next subsection.
}

\subsection{The model with change-points} Briefly, we consider the
multiple linear regression model with $T$ observations and
$m$ unknown breaks points $T_1,\ldots ,T_m$ with $1<T_1<\cdots<T_m<T$.
Here, it is important to stress that the number of change-points $m$
is known. For convenience, let $T_{0}=1$ and $T_{m+1}=T$. Namely,
let
%
\begin{equation}
\label{model} Y=\bar{Z}\delta+u,
\end{equation}
where $Y=(y_1,\ldots,y_T)^{\prime}$ is a vector of $T$ dependent
variables, $\bar{Z}$ is a $T\times(m+1)q$-matrix {of} regressors
given by $\bar{Z}=\diag(Z_1,\ldots,Z_{m+1})$ with
$Z_1=(z_{1},\ldots,z_{T_1})^{\prime}$, and for $j=2,3,\dots,m+1$,
$\mathbf{Z}_{j}= (\mathbf{z}_{T_{j-1}+1},\ldots,\mathbf
{z}_{T_{j}} )'$,
$\mathbf{z}_{T_{i-1}+1}$ {is a} $q$-column vector for
$i=1,2,\dots,m+1$. Here, $u=(u_1,\ldots,u_T)'$ is the set of
disturbances and $\delta$ is the $(m+1)q$ vector of coefficients.
Also, let $R$ be a known $k\times(m+1)q$-matrix with rank $k$,
$k\leqslant(m+1)q$ and let $r$ be a known $k$-column vector. We
consider the case where $\delta$ may satisfy or not the following
restrictions
%
\begin{equation}
\label{r1} R\delta=r.
\end{equation}
Let $\{T_1^0,\ldots,T_m^0\}$ be the true values of the break times
$\{T_1,\ldots,T_m\}$, and
$\bar{Z}^0=\diag(Z_1^0,\ldots,\allowbreak Z_{m+1}^0)$, where
$Z_i^0=(z_{T_{i-1}^0+1},\ldots,z_{T_i^0})'$. Set
$\delta= (\delta'_{1},\delta'_{2},\dots,\delta'_{m+1} )'$
where for $i=1,2,\dots,m+1$ $\delta_{i}$ is a $q$-column vector.

To estimate the unknown parameters
$(\delta'_1,\ldots ,\delta'_{m+1},T_1,\ldots,T_{m+1})'$ based only on the
sample information given in $\{Y,Z\}$, one can use the least squares
principle as described, for example in Perron and Qu \cite
{pq2006}. Also,
in case the restriction in (\ref{r1}) holds, it is common to use the
restricted least squares methods in order to estimate the target
parameter. This gives the restricted estimator (RE) of
$(\delta,T_1,\ldots,T_m)$. In particular, concerning the change-points,
let $\{\tilde{T}_1,\ldots,\tilde{T}_m\}$ denote the RE of the true
change points from restricted OLS and let
$\{\hat{T}_1,\ldots,\hat{T}_m\}$ be the unrestricted estimators (UE).
Also, let $\hat{\delta}$ and $\tilde{\delta}$ be, respectively, the UE
and RE for the regression coefficients $\delta$. Then, following the
framework in Perron and Qu \cite{pq2006}, let
$\SSR_T^R(T_1,\ldots,T_m)$ and
$\SSR_T^U(T_1,\ldots,T_m)$ be the sum of square residuals from the RE
and UE OLS regression evaluated at the partition $\{T_1,\ldots,T_m\}$,
respectively. We have
%
\begin{eqnarray}
\label{tssr} (\tilde{T}_1 ,\ldots,\tilde{T}_m)&=&\arg\min_{T_1,\ldots, T_m}\SSR_T^R(T_1,
\ldots,T_m), \nonumber\\[-8pt]\\[-8pt]
(\hat{T}_1 ,\ldots,\hat{T}_m)&=&
\arg\min_{T_1,\ldots, T_m}\SSR_T^U(T_1,
\ldots,T_m).\nonumber
\end{eqnarray}
The optimality of the proposed method is based on the asymptotic
properties of the UE and RE. In particular, in
Section~\ref{sec:jointnorm}, we establish as a preliminary step the
joint asymptotic normality of the UE and RE. To this end, we present
below the regularities conditions. To simplify the notation, let the
$\mathcal{L}_2$-norm of random matrix $X$ be defined by $\| X \|_2=
({\sum_a\sum_b} \E|X_{a,b}|^2)^{1/2}$, and let $\{{
\mathcal{F}}_{i}, i=1,2,\dots\}$ be a filtration. Also, let
$\mathrm{o}_{p}(a)$ denote a random quantity such that $\mathrm{o}_{p}(a)/a$ converges
in probability to 0, let $\mathrm{O}_{p}(a)$ denote a random quantity such
that $\mathrm{O}_{p}(a)/a$ is bounded in probability. Similarly, let $\mathrm{o}(a)$
denote a nonrandom quantity such that $\mathrm{o}(a)/a$ converges to 0, let
$\mathrm{O}(a)$ denote a nonrandom quantity such that $\mathrm{O}(a)/a$ is bounded. We
also use the notations $\displaystyle \mathop{\mathop{\xrightarrow}_{T\rightarrow
\infty}}^{\mathrm{d}}$ and $\displaystyle \mathop{\mathop{\xrightarrow}_{T\rightarrow
\infty}}^{\mathrm{P}}$ to stand for convergence in distribution and
convergence in probability respectively.

\begin{Assumptions*}[(Regularity conditions)]
\begin{enumerate}[($\mathcal{A}_3$)]
\item[($\mathcal{A}_1$)] Let $L_p=(T_{p+1}^0-T_p^0)$,
$p=1,\ldots,m$, then
$(1/L_p)\sum_{t=T_p^0+1}^{T_p^0+[L_pv]}z_tz_t^{\prime}\stackrel
{p}{\longrightarrow}
Q_p(v)$ a nonrandom positive definite matrix uniformly in $v\in
[0,1]$. Besides, there exists an $L_0>0$ such that for all ${
L_{p}}>L_0$, the minimum eigenvalues of $(1/{
L_{p}})\sum_{t=T_p^0+1}^{T_p^0+{ L_{p}}}z_tz_t^{\prime}$ and of
$(1/{ L_{p}})\sum_{t=T_p^0-{ L_{p}}}^{T_p^0}z_tz_t^{\prime}$ are
bounded away from 0.

\item[($\mathcal{A}_2$)] The matrix
$\sum_{t=i_1}^{i_2}z_tz_t^{\prime}$ is invertible for $0\leqslant
i_2-i_1\leqslant\epsilon_{0} T$ for some $\epsilon_{0}>0$.

\item[($\mathcal{A}_3$)] $T_p^0=[T\lambda_p^0]$, where
$p=1,\ldots,m+1$ and $0<\lambda_1^0<\cdots<\lambda_m^0<\lambda_{m+1}^0=1$.\vspace*{2pt}

\item[($\mathcal{A}_4$)] The minimization problem defined by
(\ref{tssr}) is taken over all possible partitions such that
$T_i-T_{i-1}>\tau T$ $(i=1,\ldots,m+1)$ for some $\tau>0$.

\item[($\mathcal{A}_5$)] For each { segment},
$(T_{p-1}^0,T_p^0)$, $p=1,\ldots,m+1$, set ${X_{pi}=
T^{-1/2}z_{T_{p-1}^0+i}u_{T_{p-1}^0+i}}$ and set ${\mathcal
{F}}_{p,i}={\mathcal{F}}_{T_{p-1}^0+i}$. We assume that
$\{X_{pi}, \mathcal{F}_{p,i}\}$ forms a $L^2$-mixingale array of size
$-1/2$. That is, there exist nonnegative constants
$\{c_{pi}:i\geqslant1\}$ and $\psi(j)$, $j\geqslant0$ such that
$\psi(j)\downarrow0$ as $j\rightarrow\infty$ and for $i\geqslant
1$, $j\geqslant0$, with
\begin{eqnarray*}\label{a6-3}
\bigl\| \E(X_{pi}|\mathcal{F}_{p,i-j})\bigr\|_2&\leqslant&
c_{pi}\psi(j),\\
 \bigl\| X_{pi}-\E(X_{pi}|
\mathcal{F}_{p,i+j})\bigr\|_2&\leqslant& c_{pi}\psi(j+1),\qquad
\psi(j)=\mathrm{O}\bigl(j^{-1/2-\epsilon}\bigr)
\end{eqnarray*}
for some $\epsilon>0$. Also, let $L_p={ T_{p+1}^{0}}-{ T_{p}^{0}}$,
and define $l_p$, $b_p$ and $r_p=[L_p/b_p]$\vspace*{2pt} such that $b_p\geqslant
l_p+1$, $l_p\geqslant1$, $b_p\leqslant L_p$. We assume that as
$b_p\displaystyle \mathop{\xrightarrow}_{L_p\rightarrow\infty}\infty$,
$l_p\displaystyle \mathop{\xrightarrow}_{L_p\rightarrow\infty} \infty$,
$b_p/L_p\rightarrow0$, and $l_p/b_p\rightarrow0$.

\item[($\mathcal{A}_6$)] For $p=1,\ldots,m+1$, for $s=1,\ldots,q$,
$\{X_{pi,s}^2/c_{pi}^2, i=1, 2, \ldots\}$ is uniformly integrable;
\[
\max_{1\leqslant i\le
L_p}c_{pi}=\mathrm{o} \bigl(b_p^{-1/2} \bigr);\qquad
{\sum_{i=1}^{r_p} \Bigl(\max_{(i-1)b_p+1\leqslant t\leqslant
ib_p}c_{pt} \Bigr)^2
=\mathrm{O}\bigl(b_p^{-1}\bigr)}
\]
 and
 \[
 {\sum_{i=1}^{r_p} \Biggl(\sum_{t=(i-1)b_p+l_p+1}^{ib_p}X_{pt} \Biggr)
 \Biggl(\sum_{t=(i-1)b_p+l_p+1}^{ib_p}X_{pt} \Biggr)' \displaystyle \mathop{\mathop{\xrightarrow}
_{L_p\rightarrow\infty}}^{p}\Sigma_p}.
\]
\end{enumerate}

Moreover, let
$V_{j,i}=\sum_{t=(i-1)b_{j}+l_{j}+1}^{ib_{j}}X_{j,t}$,
$j=1,2,\dots,m+1$. Let $r_{(1)}={\min_{1\leqslant j\leqslant
m}}(r_{p_j})$, let $r_{(m)}={\max_{1\leqslant j\leqslant
m}}(r_{p_j})$, and let $L_{\mathrm{min}}=\min(L_1,\ldots,L_{m+1})$. We have
\begin{enumerate}[2.]
\item${\sum_{i=r_{(1)}+1}^{r_{(m)}}
 (\max_{(i-1)b_{j}+1\leqslant t\leqslant ib_{j}}c_{jt}
)^2=\mathrm{o}(b_{j}^{-1})}$,
$j=1,2,\dots,m+1$.
\item${\sum_{i=1}^{r_{(1)}}}(V_{1,i}',V_{2,i}',\dots,V_{m+1,i}')'
(V_{1,i}',
V_{2,i}',\dots, V_{m+1,i}')\displaystyle \mathop{\mathop{\xrightarrow}_{L_{\mathrm{min}}\rightarrow
\infty}}^{p}\Omega$,
where $\Omega$ is nonrandom positive definite matrix.
\end{enumerate}
\end{Assumptions*}

For the interpretation of Assumptions
$(\mathcal{A}_1)$--$(\mathcal{A}_{4})$, we refer to Perron and Qu \cite{pq2006}. In summary, Assumptions $(\mathcal{A}_{1})$ and
$(\mathcal{A}_{2})$ are usually imposed in multiple linear
regressions with structural changes. Further,
Assumption $(\mathcal{A}_{3})$ guarantees to have asymptotically
distinct change points and Assumption $(\mathcal{A}_{4})$ { puts a
lower bound on the distance between breaks}. { As mentioned in
Perron and Qu \cite{pq2006}, this assumption is stronger than
the similar
condition literature. As justified in the quoted paper, this is the
cost needed to allow the heterogeneity and serial correlation in the
errors.} Assumptions $(\mathcal{A}_{5})$--$(\mathcal{A}_{6})$ are
needed to establish the asymptotic normality of the UE. Note that
Assumption $(\mathcal{A}_{5})$ considers the case of mixingale
random variables, which allow both the regressors and the errors in
each break to be a form of different distributions and
asymptotically weak dependencies.

\section{The joint asymptotic distribution of the UE and RE}
\label{sec:jointnorm} In this section, we derive the asymptotic
joint normality for the restricted and unrestricted OLS. Under
Assumptions $(\mathcal{A}_1)$--$(\mathcal{A}_{4})$,
$T^{-1} \bar{Z}^{0\prime}\bar{Z}^0$ converges in probability to a
nonrandom $q(m+1)\times q(m+1)$-positive and definite matrix.
Hereafter, we denote this matrix by $\Gamma$. Also, under
Assumption $(\mathcal{A}_{6})$,
$T^{-1} (\bar{Z}^{0\prime}uu^{\prime}\bar{Z}^0)$ converges in
probability to $\Omega$, which is a nonrandom $q(m+1)\times
q(m+1)$-positive and definite matrix. Further, under Assumptions
$(\mathcal{A}_{5})$--$(\mathcal{A}_{6})$, we establish the
following lemma which is crucial in establishing the joint
asymptotic of the UE and RE.

%
\begin{lea}\label{lemma01}
Under Assumptions $(\mathcal{A}_1)$--$(\mathcal{A}_{6})$,
$T^{-1/2}\bar{Z}^{0\prime}u\displaystyle \mathop{\mathop{\xrightarrow}_{T\rightarrow
\infty}}^{d}\mathcal{N}_{(m+1)q}(0,\Omega)$.
\end{lea}

The proof is given in the Appendix \ref{sec:appendols}. Also, note that if the
restriction in (\ref{r1}) does\vspace*{2pt} not hold, the asymptotic distribution
of $\tilde{\delta}$ may { degenerate}. Thus, in order to derive the
joint asymptotic normality, we consider the following sequence of
local alternative,
%
\begin{equation}
\label{h1t} H_{1T}: R\delta=r+\frac{\mu}{\sqrt{T}},\qquad  T=1,2,\ldots,
\end{equation}
with $\|\mu\|<\infty$. To simplify the notation, let $\hat{\delta}$
and $\tilde{\delta}$ denote, respectively, the UE and RE of $\delta$.
Let $J_0=\Gamma^{-1}R^{\prime}(R\Gamma^{-1}R^{\prime})^{-1}$, and
let $I_m$ denote $m\times m$ identity matrix. Further, let
\begin{eqnarray*}
\mu_1&=&-J_0\mu,\qquad  \Sigma_{11}=
\Gamma^{-1}\Omega\Gamma^{-1},\qquad  \Sigma_{12}=
\Gamma^{-1}\Omega\Gamma^{-1}\bigl(I_{(m+1)q}-R^{\prime
}J_0^{\prime}
\bigr),
\\
\Sigma_{21}&=&\Sigma_{12}^{\prime},\qquad
\Sigma_{22}= (I_{(m+1)q}-J_0R)\Gamma^{-1}
\Omega\Gamma^{-1}\bigl(I_{(m+1)q}-R^{\prime
}J_0^{\prime}
\bigr),
\\
\Lambda_{11}&=&J_0R\Sigma_{11}R^{\prime}J_0^{\prime},\qquad
\Lambda_{12}=J_0R\Sigma_{12}, \qquad 
\Lambda_{21}=\Lambda_{12}^{\prime},\qquad
\Lambda_{22}=\Sigma_{22}.
\end{eqnarray*}

%
\begin{lea}\label{norm1}
Under Assumptions $(\mathcal{A}_{1})$--$(\mathcal{A}_{6})$, and the
sequence of local alternative in (\ref{h1t}),
\begin{eqnarray*}
\lleft( %
\begin{array} {c} \sqrt{T}\bigl(\hat{\delta}-
\delta^0\bigr)
\\
\sqrt{T}\bigl(\tilde{\delta}-\delta^0\bigr)
\end{array} %
 \rright) &\displaystyle \mathop{\mathop{\longrightarrow}_{T\rightarrow\infty}}^{d}&
\lleft( %
\begin{array} {c} \epsilon_3
\\
\epsilon_4
\end{array} %
 \rright)\sim\mathcal{N}_{2(m+1)q}\lleft(
 \left( %
\begin{array} {c} 0
\\
\mu_1
\end{array} %
 \right), \left( %
\begin{array} {c@{\quad}c} \Sigma_{11}&
\Sigma_{12}
\\
\Sigma_{21}&\Sigma_{22}
\end{array} %
 \right)
 \rright);
\\
\lleft( %
\begin{array} {c} \sqrt{T}(\hat{\delta}-\tilde{\delta})
\\
\sqrt{T}\bigl(\tilde{\delta}-\delta^0\bigr)
\end{array} %
 \rright) &\displaystyle \mathop{\mathop{\longrightarrow}_{T\rightarrow\infty}}^{d}&
\lleft( %
\begin{array} {c} \epsilon_5
\\
\epsilon_4
\end{array} %
 \rright)\sim\mathcal{N}_{2(m+1)q}\lleft(
 \left( %
\begin{array} {c} -\mu_1
\\
\mu_1
\end{array} %
\right), \left( %
\begin{array} {c@{\quad}c} \Lambda_{11}&
\Lambda_{12}
\\
\Lambda_{21}&\Lambda_{22}
\end{array} %
 \right)
 \rright).
\end{eqnarray*}
%
\end{lea}

From the above result, it should be noted that
$(\epsilon_5,\epsilon_4)'$, the limit in distribution $(\sqrt
{T}(\hat{\delta}-\tilde{\delta}),\sqrt{T}(\tilde{\delta}-\delta^0))$
are not uncorrelated as for example in Saleh \cite{sal}, Theorem~3,
page 375, Hossain \emph{et al.} \cite{ahm2009}, among others. { Further, note
that $\Lambda_{11}$ and $\Lambda_{22}$ are not positive definite
matrices as the case in Judge and Mittelhammer \cite{judg2004}.}
Because of
that, the construction of shrinkage-type estimators cannot be
obtained by applying the results given in the literature.

\section{Shrinkage estimator and related asymptotic
properties}\label{sec:shrink} It is well known that under the
restriction in (\ref{r1}), the RE dominates in mean square error
sense the UE. However, if the restriction in (\ref{r1}) is seriously
violated, the RE performs poorly. In some scenarios, the prior
restriction in { \eqref{r1}} is subjected to some uncertainty that
may be induced by the change in the phenomenon underlying the
regression model in (\ref{model}). Under such an uncertainty, it is
of interest to propose a statistical method which combine in optimal
way the sample information and an uncertain information given
in (\ref{r1}).

In this section, we introduce a class of shrinkage estimators which
encloses the UE, RE as well as Stein-type estimator, and positive
part Stein-type estimator. To simplify some notations,\vspace*{2pt} let
$A=R'(R\Gamma^{-1}\Omega\Gamma^{-1}R')^{-1}R$, and
$\hat{A}=R'(R\hat{\Gamma}^{-1}
\hat{\Omega}\hat{\Gamma}^{-1}R')^{-1}R$, where $\hat{\Omega}$ and
$\hat{\Gamma}$ denote consistent estimators of $\Omega$ and
$\Gamma$, respectively. Also, as in Nkurunziza \cite
{n2012a}, let $h$ be
continuous (except on a number of finite points), real-valued and
integrable function (with respect to the Gaussian measure). We
consider the following class of estimators
%
\begin{eqnarray}
\label{classofSE} { \hat{\beta}}(h)=\tilde{\delta}+h \bigl(T(\tilde{\delta}-\hat {
\delta})'\hat{A} (\tilde{\delta}-\hat{\delta}) \bigr) (\hat{\delta}-
\tilde {\delta} ).
\end{eqnarray}
It should be noted that for the case where $h\equiv0$, ${
\hat{\beta}}(0)$ is the RE $\tilde{\delta}$. Also, if $h\equiv1$,
we have the UE, that is, ${ \hat{\beta}}(1)=\hat{\delta}$. Further, by
choosing a suitable $h$ one can get the pretest estimators as given
for example in Saleh \cite{sal}, Hossain \emph{et al.}
\cite{ahm2009}, among
others. Finally, the James--Stein estimator $\hat{\delta}^s$ and
Positive-Rule Stein estimator $\hat{\delta}^{s+}$ are members of the
class in \eqref{classofSE}. Indeed,{ let $k$ denote the rank of the
matrix $\mathbf{R}$ as defined in \eqref{r1}.} By taking
$h(x)=1-(k-2)/x$, $x>0$, and
$h(x)=\max \{0, 1-(k-2)/x \}$, $x>0$ we get
$\hat{\delta}^s$ and $\hat{\delta}^{s+}$, respectively.
More precisely, we have
$\hat{\delta}^s=\tilde{\delta}+ (1-\frac{k-2}{\psi} )
(\hat{\delta}-\tilde{\delta})$,
$\hat{\delta}^{s+}=\tilde{\delta}+ (1-\frac{k-2}{\psi} )^{+}
(\hat{\delta}-\tilde{\delta})$
where
$\psi=T(\tilde{\delta}-\hat{\delta})'\hat{A}(\tilde{\delta
}-\hat{\delta})$,
with $x^+=\max(0,x)$.

In order to evaluate the performance of the proposed estimators, we
consider the quadratic loss function
$L(\theta,d)=(d-\theta)'W(d-\theta)$, where $W$ is a symmetric
nonnegative definite matrix, and use the asymptotic distributional
risk (ADR) as defined, for example, in Saleh \cite{sal}.
For the
convenience of the reader, we recall  that the ADR of an
estimator $\hat{\boldsymbol{\theta}}$ is defined as
$\ADR (\hat{\boldsymbol{\theta}},\boldsymbol{\theta
};\mathbf{W} )
=\E [
\boldsymbol{\rho}'_{0}\mathbf{W}\boldsymbol{\rho}_{0} ]$,
with $\boldsymbol{\rho}_{0}$ the limit in distribution of
$\sqrt{T} (\hat{\boldsymbol{\theta}}-\boldsymbol{\theta
} )$ as $T$
tends to infinity, and $\mathbf{W}$ is a certain weight nonnegative
definite matrix.

In the sequel, we set $\Delta=\mu_1'A\mu
_1$ and assume that
the weight matrix $W$ satisfies $W=A^{1/2}W^*A^{1/2}$, with $W^*$
a symmetric nonnegative definite matrix. We establish below a lemma
which gives the ADR of estimators which are members of the class
in \eqref{classofSE}. Briefly, the derivation of this lemma is based
on the identity, established in Appendix \ref{appC}, which generalizes
Theorem~2 in Judge and Bock \cite{judg1978}. In particular,
this lemma is
useful in deriving ADR of $\hat{\delta}$, $\tilde{\delta}$,
$\hat{\delta}^s$ and $\hat{\delta}^{s+}$.

%
\begin{lea} \label{adrclass}
Suppose that Assumptions $(\mathcal{A}_1)$--$(\mathcal{A}_{6})$ and
the sequence of local alternative in \eqref{h1t} hold. Then
%
\begin{eqnarray}
& & \ADR\bigl(\hat{\beta}(h),\delta^0,W\bigr)\nonumber\\
&&\quad =\ADR\bigl(\tilde{
\delta},\delta ^0,W\bigr) -2\E\bigl[h\bigl(\chi_{k+2}^2(
\Delta)\bigr)\bigr]\mu_1'W\mu_1
\nonumber
\\
& &\qquad {} -2\E\bigl[h\bigl(\chi_{k+2}^2(\Delta)\bigr)\bigr]
\mu_1'A\Lambda_{12}W\mu_1+2\E
\bigl[h\bigl(\chi_{k+2}^2(\Delta)\bigr)\bigr]\trace(
\Lambda_{12}W\Lambda _{11}A)
\\
& &\qquad {} +2\E\bigl[h\bigl(\chi_{k+4}^2(\Delta)\bigr)\bigr]
\mu_1'A\Lambda_{12}W\mu_1\nonumber\\
&&\qquad {} +\E
\bigl[h^2\bigl(\chi_{k+2}^2(\Delta)\bigr)\bigr]
\trace(W\Lambda_{11}) +\E\bigl[h^2\bigl(
\chi_{k+4}^2(\Delta)\bigr)\bigr]\mu_1'W
\mu_1.
\nonumber
\end{eqnarray}
\end{lea}

\begin{pf}
The proof of this lemma follows directly by combining
Lemma~\ref{norm1}, Theorem~\ref{thmb1} and Lemma~\ref{lema2}.
\end{pf}

From Lemma~\ref{adrclass}, by taking $h(x)=1$, $h(x)=0$,
$h(x)=1-\frac{k-2}{x}$ and
$h(x)=\max \{0, (1-\frac{k-2}{x}) \}$, we establish the
following corollary which gives the ADR of the estimators
$\hat{\delta}$, $\tilde{\delta}$, $\hat{\delta}^s$ and
$\hat{\delta}^{s+}$, respectively.

%
\begin{cor}\label{adrestimators}
Suppose that the conditions of Lemma~\ref{adrclass} hold, then
\begin{eqnarray}\label{adrr}
 &&\ADR\bigl(\hat{\delta},\delta^0,W\bigr)\nonumber \\
 &&\quad =\trace\bigl(W
\Gamma^{-1}\Omega\Gamma^{-1}\bigr),\nonumber
\\
 && \ADR\bigl(\tilde{\delta},\delta^0,W\bigr)\nonumber\\
 &&\quad  =
\trace \bigl[W(I_{q(m+1)}-J_0R)\Gamma^{-1}
\Omega\Gamma ^{-1}\bigl(I_{q(m+1)} -R^{\prime}J_0^{\prime}
\bigr) \bigr] +\mu_1^{\prime}W\mu_1,\nonumber
\\
&&\ADR\bigl(\hat{\delta}^s,\delta^0,W\bigr)\nonumber\\
&&\quad =\ADR\bigl(
\hat{\delta},\delta^0,W\bigr)- 2(k-2)\E\bigl[\chi_{k+2}^{-2}(
\Delta)\bigr]\trace\bigl(W(\Lambda_{11}+\Lambda _{12})\bigr)\nonumber
\\
& &\qquad  { } +\bigl(k^2-4\bigr)\E\bigl[\chi_{k+4}^{-4}(
\Delta)\bigr]\mu'_1 W\mu_1+(k-2)^2
\E\bigl[\chi_{k+2}^{-4}(\Delta)\bigr]\trace(W\Lambda
_{11})\nonumber
\\
& &\qquad  { } +4(k-2)\E\bigl[\chi_{k+4}^{-4}(\Delta)\bigr]
\mu'_1 A\Lambda_{12} W\mu_1,\nonumber
\\
\label{term1}
&&\ADR\bigl(\hat{\delta}^{s+},\delta^0,W\bigr)\nonumber\\[-8pt]\\[-8pt]
&&\quad =\ADR\bigl(
\hat{\delta}^s,\delta ^0,W\bigr)
\nonumber
\\
&&\qquad {}+2\E\bigl(I\bigl(\chi^2_{k+2}(\Delta)<k-2\bigr)-(k-2)
\chi^{-2}_{k+2}(\Delta )I\bigl(\chi ^2_{k+2}(
\Delta)<k-2\bigr)\bigr) \mu_1'W\mu_1
\nonumber
\\
&&\qquad {}+2\E\bigl(I\bigl(\chi^2_{k+2}(\Delta)<k-2\bigr)-(k-2)
\chi^{-2}_{k+2}(\Delta )I\bigl(\chi ^2_{k+2}(
\Delta)<k-2\bigr)\bigr) \mu_1'A\Lambda_{12}W
\mu_1
\nonumber
\\
&&\qquad {}-2\E\bigl(I\bigl(\chi^2_{k+2}(\Delta)<k-2\bigr)-(k-2)
\chi^{-2}_{k+2}(\Delta )I\bigl(\chi ^2_{k+2}(
\Delta)<k-2\bigr)\bigr) \trace(W\Lambda_{12})
\nonumber
\\
&&\qquad {}-2\E\bigl(I\bigl(\chi^2_{k+4}(\Delta)<k-2\bigr)-(k-2)
\chi^{-2}_{k+4}(\Delta )I\bigl(\chi ^2_{k+4}(
\Delta)<k-2\bigr)\bigr) \mu_1'A\Lambda_{12}W
\mu_1
\nonumber
\\
&&\qquad {}-\E\bigl(I\bigl(\chi^2_{k+2}(\Delta)<k-2\bigr)-2(k-2)
\chi ^{-2}_{k+2}(\Delta)I\bigl(\chi^2_{k+2}(
\Delta)<k-2\bigr)
\nonumber
\\
&&\qquad \hphantom{{}-\E\bigl(}{}+(k-2)^2\chi^{-4}_{k+2}(\Delta)I\bigl(
\chi^2_{k+2}(\Delta )<k-2\bigr)\bigr)\trace(W
\Lambda_{11})
\nonumber
\\
&&\qquad {}-\E\bigl(I\bigl(\chi^2_{k+4}(\Delta)<k-2\bigr)-2(k-2)
\chi ^{-2}_{k+4}(\Delta)I\bigl(\chi^2_{k+4}(
\Delta)<k-2\bigr)
\nonumber
\\
&&\qquad\hphantom{{}-\E\bigl(} {}+(k-2)^2\chi^{-4}_{k+4}(\Delta)I\bigl(
\chi^2_{k+4}(\Delta)<k-2\bigr)\bigr)\mu _1'W
\mu_1.\nonumber
\end{eqnarray}
\end{cor}

It should be noted that the expressions in
Corollary~\ref{adrestimators} are more general than that, for
example, in Saleh \cite{sal}, page 377, and Hossain \emph{et al.} \cite{ahm2009}
for which $\Lambda_{12}=\mathbf{0}$.

From Corollary~\ref{adrestimators}, we establish the following
corollary which shows that shrinkage estimators dominate the UE. It
is noticed that, due to the asymptotic dependance structure between
the shrinking factor and the restricted estimator, the above
dominance condition looks quite complicated. To simplify the
notation, let $\Ch_{\max}(\Pi)$ denote the largest
eigenvalue of $\Pi$, and let $\Ch_{\mathrm{min}}(\Pi)$ denote the
smallest eigenvalue of $\Pi$. Further, let
$\Pi_{0}=A^{1/2}(\Lambda_{11}+4\Lambda_{12}/(k+2))W\Lambda_{11}A^{1/2}$,
$\Pi^*=(\Pi_{0}+\Pi'_{0})/2$.

%
\begin{cor}\label{adrcompar}
Suppose that Assumptions
$(\mathcal{A}_1)$--$(\mathcal{A}_{6})$ hold, and let $W$ be
nonnegative definite matrix such that
$\trace(W\Lambda_{12})\leqslant0$,
$-\Ch_{\mathrm{min}}(W\Lambda_{11})\leqslant
\Ch_{\mathrm{min}}(W\Lambda_{12})$ and
$\trace(W(\Lambda_{11}+\Lambda_{12}))
\geqslant
\max(-\trace(W\Lambda_{12}),(k+2)\Ch_{\max}(\Pi^*)/4)$.
Then,
%
\begin{eqnarray}
\label{adre00} \ADR\bigl(\hat{\delta}^{s+},\delta^0,W
\bigr)&\leqslant&\ADR\bigl(\hat{\delta }^s,\delta^0,W
\bigr)\leqslant\ADR\bigl(\hat{\delta},\delta^0,W\bigr),\qquad  \mbox{for all
$\Delta\geqslant0$.}
\end{eqnarray}
\end{cor}

%
\begin{rem} It should be
noted that the conditions for the shrinkage estimators to dominate
the unrestricted estimator are more general than given for example
in Hossain \emph{et al.} \cite{ahm2009}, Corollary~4.2, Saleh \cite
{sal}, pages 358,
 360, 382, the relations (7.4.8), (7.4.31) and (7.8.35).

Indeed, in the quoted work, we have $\Lambda_{12}=\mathbf{0}$. In this
special case, the above condition can be rewritten as
$ \{W: \frac{\trace(W\Lambda_{11})}{\Ch_{\max
}(W\Lambda_{11})}
\geqslant
\frac{k+2}{4} \}$ and this set
contains $ \{W: \frac{\trace(W\Lambda_{11})}{
\Ch_{\max}(W\Lambda_{11})}
\geqslant
\frac{k+2}{2} \}$\vspace*{2pt} which given in
the above quoted works.
\end{rem}

%
\section{Illustrative data set and numerical evaluation}\label{sec:numerical}
\subsection{Simulation study}\label{sec:sim}

In this section, we
present some Monte Carlo simulation results to evaluate the
performances of the proposed estimators. This is done by comparing
the relative mean square efficiencies (RMSE) of the estimators with
respect to the UE, $\hat{\delta}$. Recall that
$\RMSE(\delta^*)=\risk(\hat{\delta})/\risk(\delta^*)$,
where $\delta^*$ is the proposed estimator. Note that, a relative
efficiency greater than one indicates the degree of superiority of
the proposed estimator over $\hat{\delta}$. To save the space of
this paper, we report only two cases.

\textit{Case} 1: the number of unknown parameters is
small, with $m=3$, $q=2$;
$\delta^{0}= (\delta_1^{0'},\delta_2^{0'},\delta_3^{0'},\delta
_4^{0'} )'$
with $\delta_1^0=\delta_3^0= (1,2 )'$ and
$\delta_2^0=\delta_4^0=\mathbf{0}$ (i.e., the zero vector), and the
sample sizes are set to be $T=40$ with the change points given by
$(10,20,30,40)$. Also, we set $T=100$ with the change-points
$(25,50,75,100)$. Further, the restriction is such that
$\mathbf{R}=
[E_{1},E_{2},E_{3},E_{4},-E_{1},-E_{2},\allowbreak E_{5},E_{6} ]$
where, for $j=1,2,\dots,6$, $E_{j}$ is a 6-column vector with all
components equal to zero except the $j$th component which
equal to 1.

\textit{Case} 2: the number of unknown parameters is
relative large by setting $m=4$, $q=5$,
$\delta^{0}= (\delta_1^{0'},\delta_2^{0'},\delta_3^{0'},\delta_4^{0'},
\delta_5^{0'} )'$ with
$\delta_1^0=\delta_3^0=\delta_5^0= (1,2,3,4,5 )'$,
$\delta_2^0=\delta_4^0=\mathbf{0}$ and the sample sizes are $T=100$ and
$T=500$ with the change-points $(20,40,60,80,100)$ and
$(100,200,300,400,500)$, respectively. Further, the restriction $R$
is set to be a $8\times25$ matrix with
\begin{eqnarray*}
R_{1,1}&=&R_{2,2}=R_{3,3}=R_{4,4}=R_{5,5}=R_{6,6}=R_{7,19}=R_{8,20}=1,
\\
R_{1,11}&=&R_{2,12}=R_{3,13}=R_{4,14}=R_{5,15}=-1,
\end{eqnarray*}
and the rest elements of $R$ are set to be 0.

In each case, we let $z_{T_{i}}\sim\mathcal{N}_{q}(1,\Sigma)$,
where $\Sigma$ is a $q\times q$ symmetric matrix such that
$\Sigma_{a,b}=|0.5|^{|a-b|}$. Also, we let $u_i\sim\mathcal{N}(0,
\sigma^2)$, $1\leqslant\sigma^2\leqslant2$, and compute the
related RMSE based on the 1000 replications.
%
\begin{figure}

\includegraphics{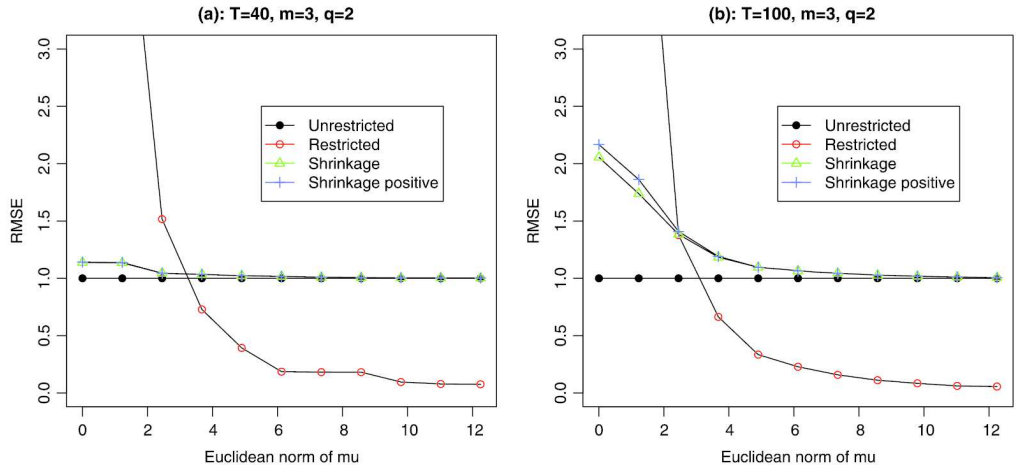}

\caption{RMSE of the restricted and shrinkage estimators (case 1).}\label{fig1}
\end{figure}
%
\begin{figure}[b]

\includegraphics{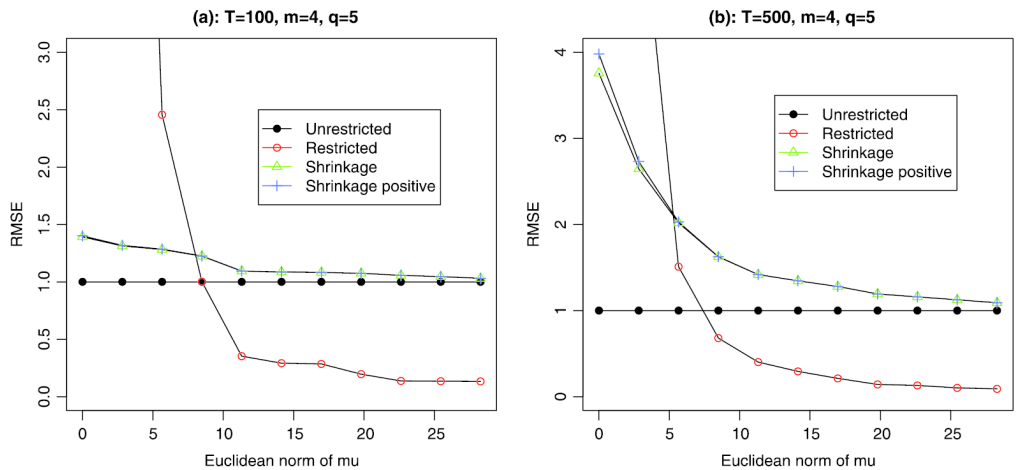}

\caption{RMSE of the restricted and shrinkage estimators (case 2).}\label{fig2}
\end{figure}

The results of the simulation studies are given in Figures~\ref{fig1} and \ref{fig2}.
In summary, the results corroborate the theoretical finding (given in Corollary~\ref{adrcompar}) for which the proposed shrinkage
estimators dominate the unrestricted estimator. {We also construct,
and present in Appendix \ref{appC}, Figures~\ref{fig3}--\ref{fig6} which give some histograms
of the UE and RE of the change points. The results given in
Figures~\ref{fig3}--\ref{fig6} suggest that both the unrestricted and the restricted
methods work well in estimating the change points.}

\renewcommand{\thefigure}{\arabic{figure}}
\setcounter{figure}{2}
\begin{figure}

\includegraphics{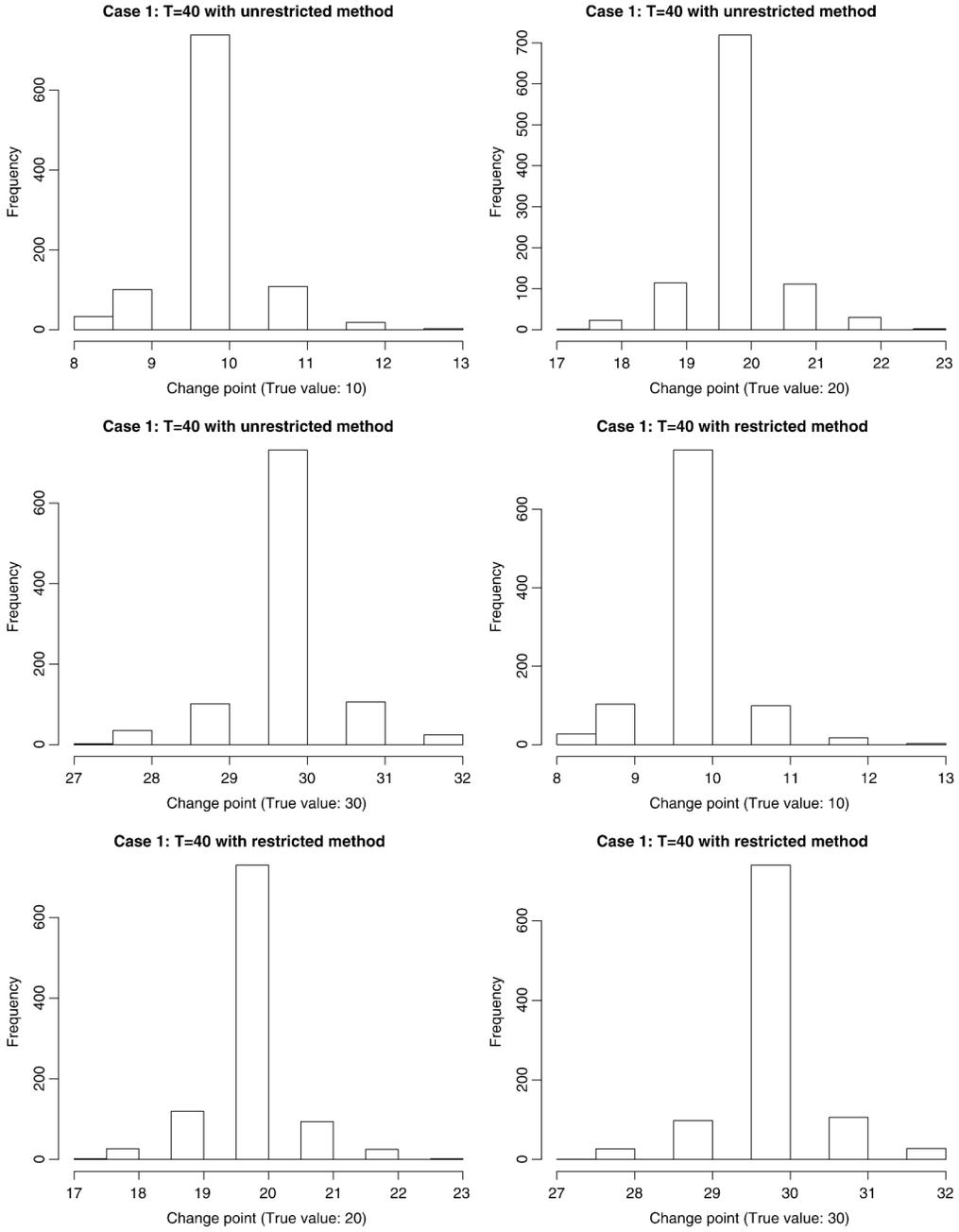}

\caption{Histograms of the UE and RE of change points (case 1 with
$T=40$).}\label{fig3}
\end{figure}
%
\begin{figure}

\includegraphics{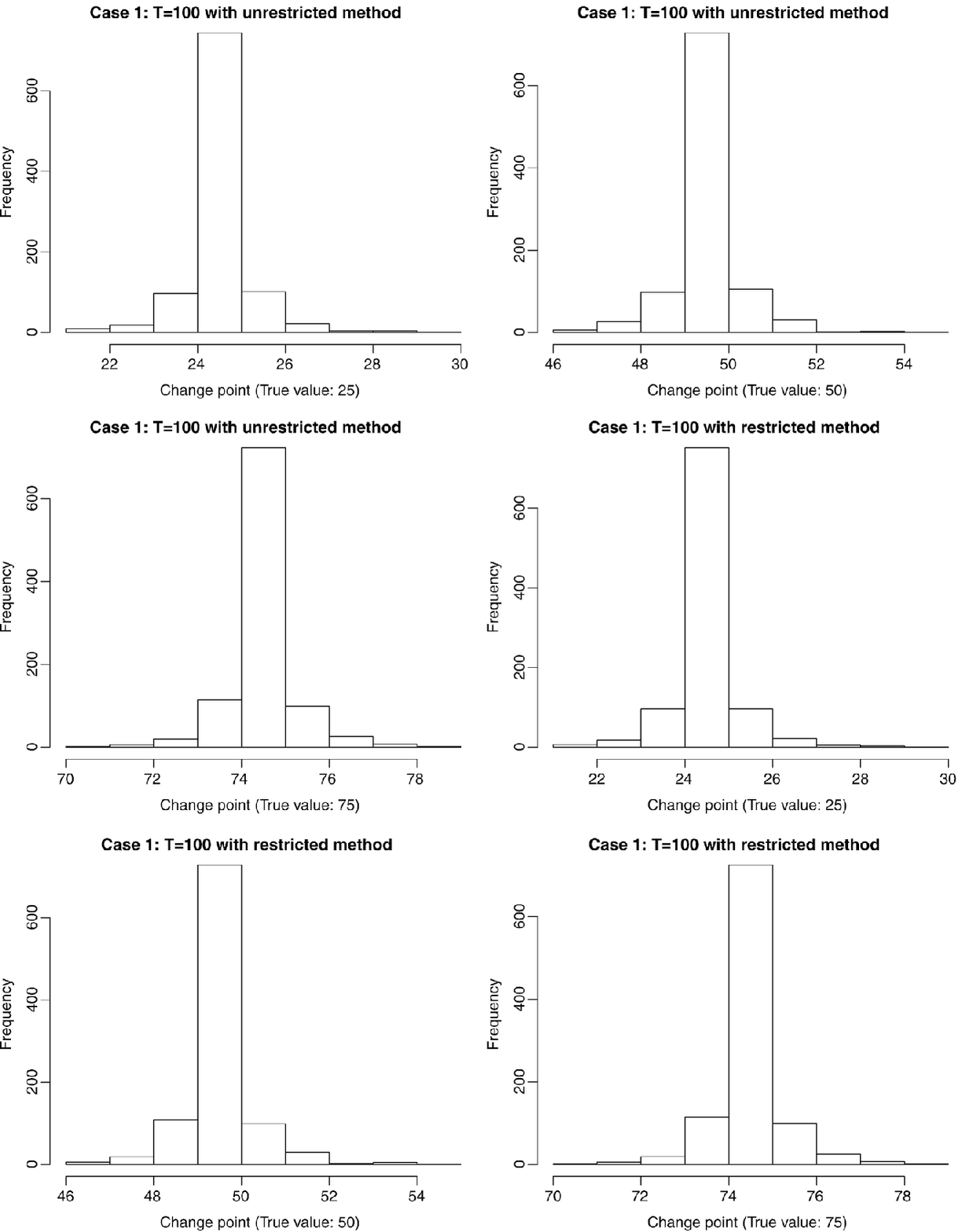}

\caption{Histograms of the UE and RE of change points (case 1 with
$T=100$).}\label{fig4}
\end{figure}
%
\begin{figure}

\includegraphics{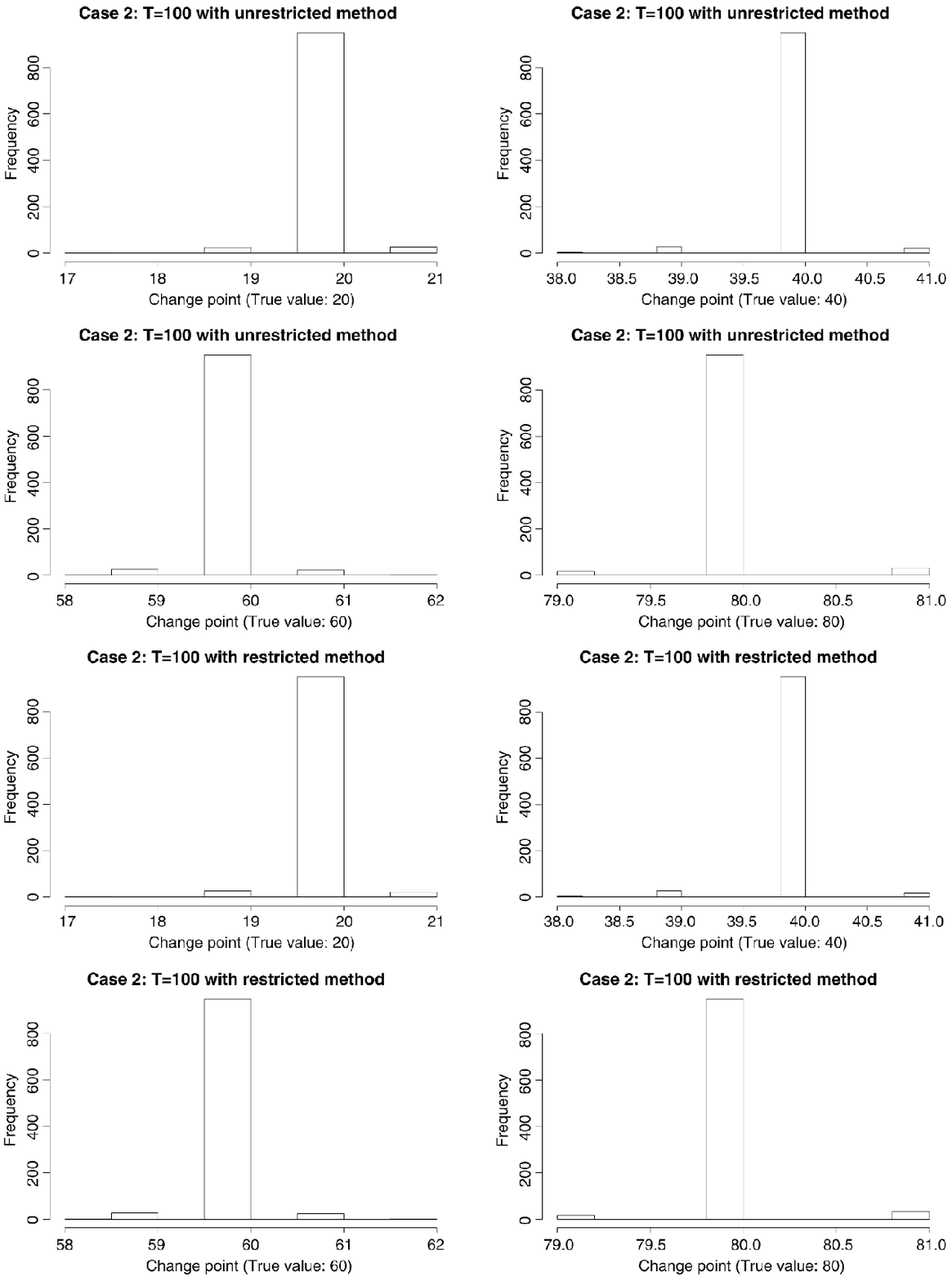}

\caption{Histograms of the UE and RE of change points (case 2 with
$T=100$).}\label{fig5}
\end{figure}
%
\begin{figure}

\includegraphics{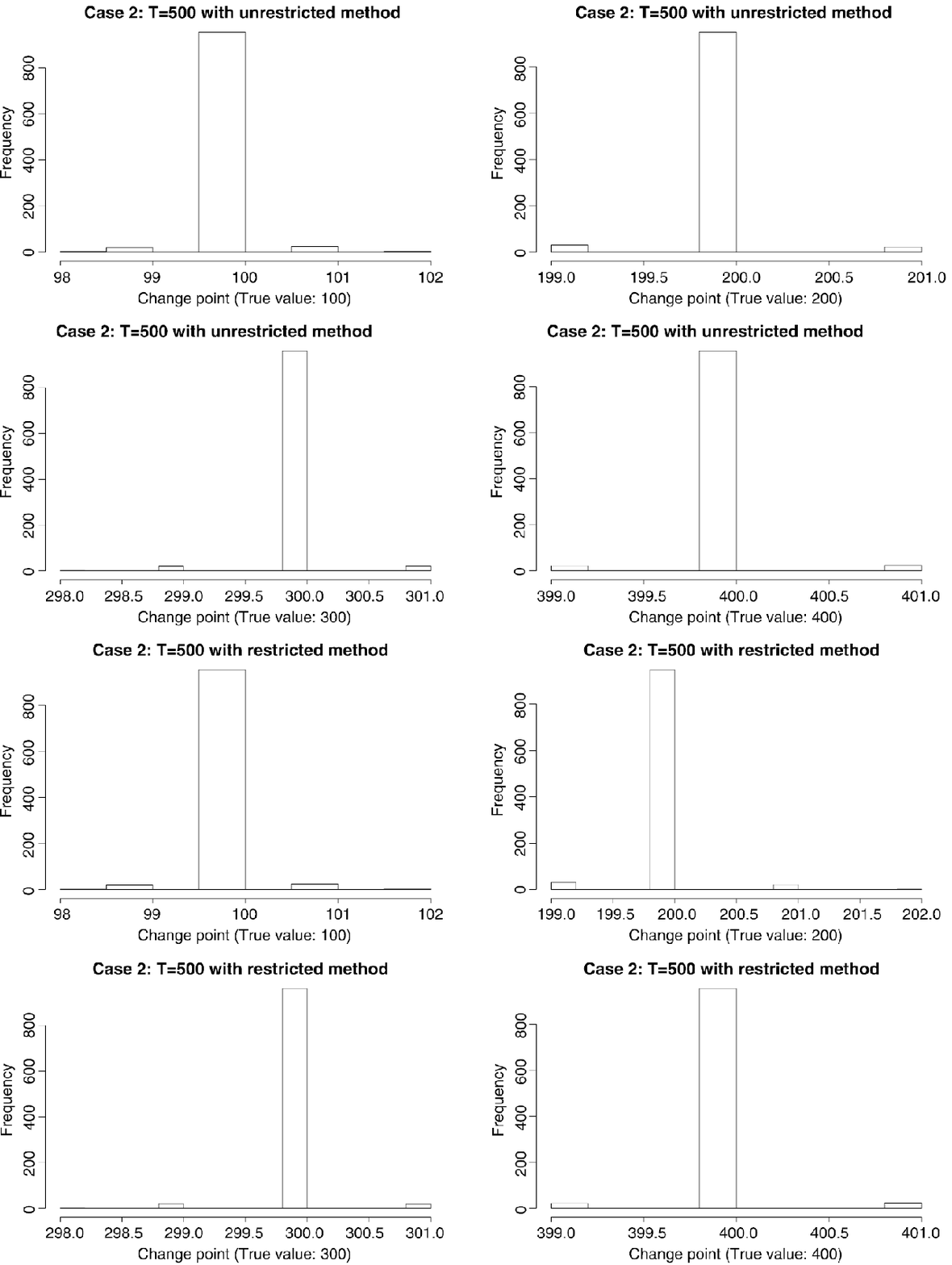}

\caption{Histograms of the UE and RE of change points (case 2 with
$T=500$).}\label{fig6}
\end{figure}

\begin{table}[b]
\tablewidth=\textwidth
\tabcolsep=0pt
\caption{Change-points and MSE}\label{table1}
\begin{tabular*}{\textwidth}{@{\extracolsep{\fill}}llld{2.8}d{1.8}d{2.8}d{2.8}@{}}
\hline
 & \multicolumn{2}{l}{Change-points} & \multicolumn
{4}{l}{MSE}\\
[-5pt]
 & \multicolumn{2}{l}{\hrulefill} & \multicolumn
{4}{l}{\hrulefill}\\
Country&(UE)& (RE) &\multicolumn{1}{l}{$\hat{\delta}$} &\multicolumn{1}{l}{$\tilde{\delta}$} &
\multicolumn{1}{l}{$\hat{\delta}^s$} &
\multicolumn{1}{l}{$\hat{\delta}^{s+}$}\\
\hline
Australia& 1907 & 1929 & 1.67004021 & 0.03936242 & 1.64839567 &
1.64839567\\
Canada &1931 & 1930 & 2.96623326 & 0.05474518 & 2.87279365 &
2.87279365\\
Denmark&1939& 1939& 3.99038175 & 0.04765026 & 3.93532691 & 3.93532691 \\
France&1943& 1943 & 12.1123258 & 0.1253509 & 11.9030741 & 11.9030741\\
Germany&1945& 1954 & 11.4218637 & 0.1704905 & 11.3279191 & 11.3279191\\
Italy&1943& 1943 & 10.2462836 & 0.1211837 & 10.2079175 & 10.2079175\\
Norway&1944& 1948 & 7.09593981 & 0.03606614 & 6.92396377 & 6.92396377\\
Sweden&1924 & 1916 & 0.72605495 & 0.02192452 & 0.70854206 & 0.70854206\\
U.K. &1918& 1919 & 0.61037392 & 0.01496282 & 0.58916536 & 0.57701458\\
U.S. &1940& 1929& 3.97869572 & 0.05967521 & 3.91443168 & 3.91443168\\
\hline
\end{tabular*}
\end{table}

\subsection{Data analysis}\label{sec:data} In this subsection, we
illustrate the application of the proposed estimation strategy to
the real data set. As a real data set, we consider a historical
($\log$) gross domestic product (GDP) data set from 1870 to 1986 for
10 different countries. This data set is used for example in
Perron and Yabu \cite{by2009}, and these authors pointed out that
most GDP series
presented in the given data set are characterized by at least one
major shift and therefore change-point model is applicable. For each
GDP series, we consider the following model:
\[
Y_t=\lleft\{ %
\begin{array} {l@{\qquad}l} \delta'_{1}
\bigl(1,t, t^{1.5},t^2 \bigr)',& \mbox{if } t=1,
\ldots,T_{1},
\\
\delta'_{2} \bigl(1,t, t^{1.5},t^2
\bigr)', &\mbox{if } t=T_{1}+1,\ldots,117, \end{array}
\rright.
\]
with $1\leqslant T_{1}\leqslant117$, for $i=1, 2$, $\delta_{i}$ is
a $4$-column vector. The uncertain restriction is given by
$R\delta=r$ with
\[
R=\lleft[ %
\begin{array} {c@{\quad}c@{\quad}c@{\quad}c@{\quad}c@{\quad}c@{\quad}c@{\quad}c} 0&0&1&0&0&0&0&0
\\
0&0&0&1&0&0&0&0
\\
0&0&0&0&0&0&1&0
\\
0&0&0&0&0&0&0&1
\\
\end{array} %
 \rright],
\]
and $r=\{0,0,0,0\}'$. In practice, the hypothesized restriction
means that the $\log(\mathrm{GDP})$ is suspected to have a linear trend. For
the given data, we first use the proposed method to calculate the
unrestricted and the restricted estimators of the
change-point $\hat{T}_{1}$ and $\tilde{T}_{1}$ as well as the
estimators $\hat{\delta}$, $\tilde{\delta}$, $\hat{\delta}^s$ and
$\hat{\delta}^{s+}$. For the change-point $T_{1}$ which is a
nuisance parameter here, we do not compute the shrinkage estimators.
The obtained unrestricted and restricted estimate of the
change-point $\hat{T}_{1}$ and $\tilde{T}_{1}$ are given in
Table~\ref{table1}. In order to save the space of this paper, we do
not report here the point estimates of $\hat{\delta},
\tilde{\delta}, \hat{\delta}^s,\hat{\delta}^{s+}$, but these values
are available upon request. Further, we calculate the MSE of each
type of estimators, by applying the bootstrap method to the
residuals. Recall that, in this paper, the change-points are treated
as the nuisance parameters. Thus, the construction of the shrinkage
estimators for the change-points is beyond the scope of this paper.

As we can see from Table~\ref{table1}, the MSE of the restricted
estimator is much smaller the MSE of the other estimators. This may
indicate that the true value of the parameter vector lies in the
neighborhood of the chosen restriction. Further, the MSE of the
proposed shrinkage estimators is smaller than the MSE of the
unrestricted estimator. The obtained result is in agreement with the
above simulation study.

\section{Conclusion}\label{sec:conclu}
The goal of this research was to derive an improved estimation
strategy for the regression coefficients in multiple linear model
with unknown change-points under uncertain restrictions. In summary,
we introduced a class of estimators which includes the UE
$\hat{\delta}$, RE $\tilde{\delta}$, James--Stein Estimator
$\hat{\delta}^s$ and Positive-Rule Stein Estimator
$\hat{\delta}^{s+}$. The main difficulty consists in the fact that
the random quantities $\tilde{\delta}-\delta$ and
$\hat{\delta}-\tilde{\delta}$ are not asymptotically uncorrelated as
this is the common case in literature. To tackle this difficulty,
we generalized (in the Appendix \ref{appC}) Theorems 1--2 in Judge and Bock \cite{judg1978}. Under the conditions more general than
that in
literature, we established that $\hat{\delta}^s$ and
$\hat{\delta}^{s+}$ dominate UE. The performance of SEs over the UE
is confirmed by the simulation studies. They also show that SEs
perform better than the RE when one moves far away from the
hypothesized restriction. It should be noticed that, in this paper,
the tools used for studying shrinkage estimators are based on
noncentral chi-squares. One of the referees suggested to
investigate if the obtained results can be improved by using more
recent tools such as Stein's unbiased risk estimate. Research on
this interesting idea is ongoing.

Another highlight of this paper consists in the fact that, in
deriving the joint asymptotic normality of the UE and RE, we relax
some conditions given in recent literature. In particular, we
considered here the condition of $L_2$-mixingale with size $-1/2$,
which allow both the regressors and the errors in each break to be a
form of different distributions and asymptotically weak
dependencies.

\begin{appendix}
\section*{Appendix}\label{app}
In this section, we give some technical proofs underlying the
results established in this paper. To set up additional notations,
let $\| A \|$ denote the Euclidean norm for vector $A$. For a matrix
$B$, let $\| B \|$ be the vector induced norm (i.e., $\| B
\|={\sup_{x\neq0}}\| Bx \|/\| x \|$).
\setcounter{section}{0}
\section{Technical results underlying the asymptotic
properties}\label{sec:techna} First, we establish the following
proposition which plays a central role in deriving the joint
asymptotic normality between the UE and RE. For the sake of
simplicity, we set $D_{i,k^*}=X_{pi}-\mathrm{E}(X_{pi}|\mathcal
{F}_{p,i+k^*})$ and set $D_{i,k^*,s}$ be
the $s$th element in $D_{i,k^*}$.

%
\begin{prn}\label{cor101}
Suppose that Assumptions $(\mathcal{A}_{5})$ and $(\mathcal{A}_{6})$
hold. Then,
\begin{eqnarray*}
& & \mathrm{E} \Biggl({\sum_{r=1}^{L_p}}(D_{i,k^{*}-1,s}-D_{i,k^{*},s})^2
\Biggr) =\sum_{i=1}^{L_p}\mathrm{E}\bigl(
\mathrm{E}^2(X_{pi,s}|\mathcal{F}_{p,i+k^*})\bigr)-
\sum_{i=1}^{L_p}\mathrm{E}\bigl(
\mathrm{E}^2(X_{pi,s}|\mathcal {F}_{p,i+k^*-1})\bigr),
\\
& & \sum_{i=1}^{L_p}\sum
_{j=l}^{i-1}\mathrm {E}\bigl[(D_{i,k^{*}-1,s}-D_{i,k^{*},s})
(D_{j,k^{*}-1,s}-D_{j,k^{*},s})\bigr]=0
\end{eqnarray*}
and
\begin{eqnarray*}
& & \sum_{i=1}^{L_p}\bigl[\mathrm{E}
\bigl(\mathrm{E}^2(X_{pi,s}|\mathcal {F}_{p,i+k^*})
\bigr) -\mathrm{E}\bigl(\mathrm{E}^2(X_{pi,s}|
\mathcal{F}_{p,i+k^*-1})\bigr)\bigr] =\sum_{i=1}^{L_{p}}
\bigl[\mathrm{E} \bigl(D_{i,k^*-1,s}^2\bigr)-\mathrm{E}\bigl(
D_{i,k^*,s}^2\bigr)\bigr].
\end{eqnarray*}
\end{prn}

\begin{pf}
One can verify that
\[
X_{pi}={\sum_{k^*=-\infty}^{\infty}}\bigl[\E(X_{pi}|\mathcal
{F}_{p,i+k^*}) -\E(X_{pi}|\mathcal{F}_{p,i+k^*-1})\bigr]\qquad \mbox{a.s.}
\]

Further, one can verify that
\begin{eqnarray*}
&&\mathrm{E} \Biggl({\sum_{r=1}^{L_p}}(D_{i,k^{*}-1,s}-D_{i,k^{*},s})^2
\Biggr)\\
&&\quad  =\sum_{i=1}^{L_p}\E\bigl[
\E^2(X_{pi,s}|\mathcal{F}_{p,i+k^*})\bigr] +\sum
_{i=1}^{L_p}\E\bigl[\E^2(X_{pi,s}|
\mathcal{F}_{p,i+k^*-1})\bigr]
\\
&&\qquad {}-2\sum_{i=1}^{L_p}\E\bigl[
\E(X_{pi,s}|\mathcal{F}_{p,i+k^*-1}) \E\bigl(\E (X_{pi,s}|
\mathcal{F}_{p,i+k^*})|\mathcal{F}_{p,i+k^*-1}\bigr)\bigr],
\end{eqnarray*}
and then, by using the properties of the conditional expected value,
we prove the first statement. For the second statement, we have
\begin{eqnarray*}
&&\sum_{i=1}^{L_p}\sum
_{j=1}^{i-1}\E\bigl[\bigl((D_{i,k^{*}-1,s}-D_{i,k^{*},s})
(D_{j,k^{*}-1,s}-D_{j,k^{*},s})\bigr)\bigr]
\\
& & \quad = \sum_{i=1}^{L_p}\sum
_{j=1}^{i-1}\E\bigl[(D_{j,k^{*}-1,s}-D_{j,k^{*},s})
\bigl(\E\bigl((D_{i,k^{*}-1,s}-D_{i,k^{*},s})|\mathcal{F}_{p,j+k^*}
\bigr)\bigr)\bigr]=0.
\end{eqnarray*}
The third statement of the proposition follows from the similar
algebraic computations.
\end{pf}

%
\begin{lea}\label{lem1}
Let $v_{L_p}^2={\sum_{i=1}^{L_p}}c_{pi}^2$ and suppose that
Assumptions $(\mathcal{A}_{5})$ and $(\mathcal{A}_{6})$ hold. Then
\[
{\sum_{s=1}^{q}}\E \Biggl({\max_{j\leqslant
L_p}} \Biggl({\sum_{i=1}^{j}}X_{pi,s} \Biggr)^2 \Biggr) \leqslant
16v_{L_p}^2 \Biggl[\sum_{k^*=0}^{\infty} \Biggl(\sum_{i=0}^{k^*}\psi
^{-2}(i) \Biggr)^{-1/2}
 \Biggr]^2.
 \]
\end{lea}

The proof follows from Proposition~\ref{cor101} and following the
similar steps as in proof of Lemma~3.2 in Mcleish \cite{mcleish1977}. By using
this lemma, one establishes the following corollary which plays a
central role in establishing the joint asymptotic normality of UE
and RE.

%
\begin{cor}\label{lema01}
Under Assumptions $(\mathcal{A}_{5})$ and $(\mathcal{A}_{6})$, then
\[
{\sum_{s=1}^{q}}\E \Biggl[ \Biggl({\sum_{i=1}^{L_p}}X_{pi,s}
\Biggr)^2 \Biggr]
=\mathrm{O} \bigl(v_{L_p}^2 \bigr).
\]
\end{cor}

\begin{pf}
From Lemma~\ref{lem1},
%
\begin{eqnarray}
\label{l1e2} {\sum_{s=1}^{q}}\E \Biggl({
\max_{j\leqslant
L_p}} \Biggl({\sum_{i=1}^{j}}X_{pi,s}
\Biggr)^2 \Biggr) \leqslant16v_{L_p}^2 \Biggl[
\sum_{k^*=0}^{\infty} \Biggl(\sum
_{i=0}^{k^*}\psi^{-2}(i)
\Biggr)^{-1/2} \Biggr]^2,
\end{eqnarray}
and then, the proof follows directly from the fact that $
{\sum_{k^*=0}^{\infty}} ({\sum_{i=0}^{k^*}}\psi^{-2}(i) )^{-1/2}
<\infty$.
\end{pf}

\begin{cor}\label{cor201}
Let $v_i^2={\sum_{(i-1) b_p+l_p+1}^{ib_p}}c_{pt}^2$ and suppose
that Assumptions $(\mathcal{A}_{5})$ and $(\mathcal{A}_{6})$ hold. Then, $
 \{{\sum_{s=1}^{q}}{\max_{j\leqslant
ib_p}} ({\sum_{t=(i-1)b_p+l_p+1}^{j}}X_{pt,s} )^{2}  /{v_{i}^2},
i=1,\ldots,r_p, r_p\geqslant1 \} $ is uniformly integrable.
In particular, $
 \{{\sum_{s=1}^{q}} ({\sum_{t=(i-1)b_p+l_p+1}^{ib_p}}X_{pt,s} )^2  /(v_{i}^2),
i=1,\ldots,r_p, r_p\geqslant1 \}$ is uniformly integrable.
\end{cor}

\begin{pf} Let $S_{j,s}={\sum_{i=1}^j}X_{pi,s}$, $s=1,\ldots,q$.
By using the same arguments as used in proof of Lemma~3.5 in
McLeish \cite{mcleish1977}, one verifies that the set
$ \{{\max_{j\leqslant
L_p}}{\sum_{s=1}^q}\frac{S_{j,s}^2}{v_{L_p}^2};L_p\geqslant
1 \} $ is uniformly integrable. This completes the proof.
\end{pf}

Further, by using Lemma~\ref{lema01}, we establish the following
proposition which is also useful in establishing the joint
asymptotic normality of UE and RE. To simplify some notations, let
$r_{\mathrm{min}}={\min_{1\leqslant p\leqslant m+1}}(r_{p})$, and let
$L_{\mathrm{min}}={\min_{1\leqslant p\leqslant m+1}}(L_{p})$. Further, let
$\mathcal{H}_i$ be the $\sigma$-field generated by
$\{U_{ib_p},U_{ib_p-1},\dots\}$, with $U_i$ are random variables
defined on $(\Omega, \mathcal{F}, P)$ such that $\mathcal
{H}_{i-1}\subseteq\mathcal{F}_{p,i-j}$, and let
$V_{pi}={\sum_{t=(i-1)b_p+l_p+1}^{ib_p}}X_{pt}$, let
$W_{p,i}=\E(V_{pi}|\mathcal{H}_i)-\E(V_{pi}|\mathcal{H}_{i-1})$,
$p=1,2,\dots,m+1$, $i=1,2,\dots, r_{\mathrm{min}}$.

%
\begin{prn}\label{convjointW}
Suppose that Assumptions $(\mathcal{A}_{5})$ and $(\mathcal{A}_{6})$
hold. Then,
\begin{eqnarray*}
&&{\sum_{i=1}^{r_{\mathrm{min}}}}\bigl[\bigl(V_{1,i}',\dots,V_{m+1,i}'\bigr)'\bigl(V_{1,i}',\dots
,V_{m+1,i}'\bigr)\\
&&\hphantom{{\sum_{i=1}^{r_{\mathrm{min}}}}\bigl[}{}-\bigl(W_{1,i}',\dots,W_{m+1,i}'\bigr)'\bigl(W_{1,i}',\dots,W_{m+1,i}'\bigr)\bigr]
\displaystyle \mathop{\mathop{\xrightarrow}_{L_{\mathrm{min}}\rightarrow\infty}}^{p}0.
\end{eqnarray*}
\end{prn}

The proof follows from Lemma~\ref{lema01} along with some algebraic
computations.

%
\begin{prn}\label{prn304}
Suppose that the conditions Proposition~\textup{\ref{convjointW}} hold.
Then,
\[
{\sum_{i=1}^{r_{\mathrm{min}}}}\bigl(W_{1,i}',W_{2,i}',\ldots,W_{m+1,i}'\bigr)'
\bigl(W_{1,i}',W_{2,i}',\ldots,W_{m+1,i}'\bigr)  \displaystyle \mathop{\mathop{\xrightarrow}_
{L_{\mathrm{min}}\rightarrow
\infty}}^{p}\Omega
\]
and
\[
{\sum_{a=1}^{m+1}}{\sum_{i=1}^{r_{a}}}{\sum_{s=1}^{q}}\E\Biggl[(W_{a,i,s})^2
\mathbb{I}\Biggl(\sum_{s=1}^{q}W_{a,i,s}^2>\epsilon\Biggr)\Biggr]\displaystyle \mathop{\xrightarrow}_{L_{\mathrm{min}}
\rightarrow\infty} 0,\qquad  \mbox{for all }\epsilon>0.
\]
\end{prn}

\begin{pf}
By using Assumption $(\mathcal{A}_{6})$ along with
Proposition~\ref{convjointW} and  Slutsky's theorem, we
establish the first statement. For the second statement, one
verifies that, for each $a=1, 2, \dots, m+1$, $\{W_{a,i},\mathcal
{H}_i\}$ is a $L_2$-mixingale array of size $-1/2$. Then the rest of
the proof follows from Corollary~\ref{cor201}.
\end{pf}

%
\section{Asymptotic normality of the UE and RE}\label{sec:appendols}
\vspace*{-12pt}%
\begin{pf*}{Proof of Lemma~\ref{lemma01}}
Note that
\[
{T^{-1/2}\bar{Z}^{0\prime}u\equiv \Biggl(\sum_{i=1}^{L_{1}}X'_{1,i},
\ldots,\sum_{i=1}^{L_{m+1}}X'_{m+1,i} \Biggr)'},
\]
 then
%
\begin{equation}\label{decom}
T^{-1/2}\bar{\mathbf{Z}}^{0\prime}u= {\sum
_{i=1}^{r_{\mathrm{min}}}}\mathbf{W}_{i}+\boldsymbol{
\Xi}^* + \Biggl( {\sum_{i=r_{\mathrm{min}}}^{r_1}}\sum
_{t=(i-1)b_{1}+1}^{{i}b_{1}}X'_{1,i},
\ldots,{\sum_{i=r_{\mathrm{min}}}^{r_{m+1}}}\sum
_{t=(i-1)b_{m+1}+1}^{{i}b_{m+1}}X'_{m+1,i}
\Biggr)',
\end{equation}
with $r_{\mathrm{min}}=\min_{1\leqslant i\leqslant
m+1}(r_{i})$ and $\boldsymbol{\Xi}^{*}= (\Xi_{1}^{*'},
\Xi_{2}^{*'},\dots,\Xi_{m+1}^{*'} )'$, where
\begin{eqnarray*}
\Xi_{j}^*=\sum_{i=1}^{r_{\mathrm{min}}}
\Biggl(V_{ji} -W_{j,i}+\sum_{t=(i-1)b_{j}+1}^{{i-1}b_{j}+l_{j}}X_{j,i}
\Biggr)+\sum_{t=r_{j}b_{j}+1}^{L_{j}}X_{p_j,t}.
\end{eqnarray*}
Further, it should be noted that, under
Assumptions ${(\mathbf{\mathcal{A}}_{4})}$ and
${(\mathbf{\mathcal{A}}_{5})}$, $T$ tends to infinity if and only if
$L_{\mathrm{min}}={\min_{1\leqslant j\leqslant m+1}}(L_{j})$ tends to
infinity.

By using Lemma~\ref{lema01} along with some algebraic computations,
we have
%
\begin{eqnarray}\label{res t1}
\bigl(\Xi_{1}^{*'}, \Xi_{2}^{*'},
\dots,\Xi_{m+1}^{*'} \bigr)'&\displaystyle \mathop{\mathop{\xrightarrow}_
{L_{\mathrm{min}}\rightarrow \infty}}^{\mathrm{P}}& 0,\nonumber
\\[-8pt]\\[-8pt]
 \Biggl(\sum_{i=r_{\mathrm{min}}}^{r_1}\sum
_{t=(i-1)b_{1}+1}^{{i}b_{1}}X'_{1,i},\dots,
\sum_{j=r_{\mathrm{min}}}^{r_{m+1}}\sum
_{t=(j-1)b_{m+1}+1}^{{j}b_{m+1}}X'_{m+1,j}
\Biggr)' &\displaystyle \mathop{\mathop{\xrightarrow}_{L_{\mathrm{min}}\rightarrow\infty}}^{
\mathrm{P}} &0.
\nonumber
\end{eqnarray}
Therefore, the proof follows from the relations \eqref{decom} and
\eqref{res t1} along with the martingale difference sequence central
limit theorem along with Slutsky's theorem.
\end{pf*}

%
\begin{prn}\label{prn3}
Under $(\mathcal{A}_1)$--$(\mathcal{A}_{6})$, we have
$\sqrt{T}(\hat{\delta}-\delta^0)\displaystyle \mathop{\mathop{\xrightarrow}_{T\rightarrow\infty}}^
{d}\epsilon_1\sim\mathcal{N}_{q(m+1)} (0,\break
\Gamma^{-1}\Omega\Gamma^{-1} )$.
\end{prn}

The proof follows by combining Lemma~\ref{lemma01} and Slutsky's
theorem.

\begin{pf*}{Proof of Proposition~\ref{norm1}}
Let
$J=(\bar{Z}^{0\prime}\bar{Z}^0)^{-1}R^{\prime}(R(\bar{Z}^{0\prime}
\bar{Z}^0)^{-1}R^{\prime})^{-1}$, we have
\begin{eqnarray*}
\bigl(\sqrt{T}\bigl(\hat{\delta}-\delta^0\bigr)',
\sqrt{T}\bigl(\tilde{\delta }-\delta^0\bigr)'
\bigr)'\doteq \bigl(I_{(m+1)q}, I_{(m+1)q}-R'J'
\bigr)'\sqrt{T}\bigl(\hat{\delta }-\delta^0\bigr) +
\bigl(0, -\mu'J' \bigr)'.
\end{eqnarray*}
Then, the first statement follows directly from
Proposition~\ref{prn3} and Slutsky's theorem, along with some
algebraic computations. For the second statement, obviously
\[
\bigl((\hat{\delta}-\tilde{\delta})', \bigl(\tilde{\delta}-\delta
^0\bigr)' \bigr)' = \bigl((I_{q(m+1)},0)',
(-I_{q(m+1)}, I_{q(m+1)})' \bigr)' \bigl(
\bigl(\hat{\delta}-\delta^{0}\bigr)', \bigl(\tilde{
\delta}-\delta ^{0}\bigr)' \bigr)'.
\]
Then, the rest of the proof follows directly from the first
statement of the proposition along with Slutsky's theorem.
\end{pf*}

%
\section{Some results for the derivation of risk functions}\label{appC}
%
\begin{thm}\label{thma1}
Let $h$ be Borel measurable and real-valued integrable function, let
$X\sim\mathcal{N}_{p}(\mu,\Sigma)$, where $\Sigma$ is a nonnegative
definite matrix with rank $k\leqslant p$. Let $A$ be a $p\times
p$-nonnegative definite matrix with rank $k$ such that $\Sigma A$ is
an idempotent matrix, $A\Sigma A=A$; $\Sigma A \Sigma=\Sigma$; and
$\Sigma A \mu=\mu$, and let $W=A^{1/2}W^*A^{1/2}$ where $W^*$ is a
nonnegative definite matrix. Then, $\E[h(X'AX)WX]=\E[h(\chi
_{k+2}^2(\mu'A\mu))]W\mu$.
\end{thm}

\begin{pf}
Let $A^{1/2\dag}$ be the Moore--Penrose pseudoinverse of $A^{1/2}$.
By the definition of Moore--Penrose pseudo-inverse, we have
$WX=A^{1/2}W^*A^{1/2}A^{1/2\dag}A^{1/2}X=WA^{1/2\dag}A^{1/2}X$, and
then,
%
\begin{eqnarray}
\E\bigl[h\bigl(X'AX\bigr)X'WX\bigr]=\E\bigl[h
\bigl(X'AX\bigr)X'A^{1/2}A^{1/2\dag}WA^{1/2\dag}A^{1/2}X
\bigr].
\end{eqnarray}
Further, since $A^{1/2}\Sigma A^{1/2}$ is a symmetric and idempotent
matrix, there exists an\vspace*{-3pt} orthogonal matrix $G$ such that $
GA^{1/2}\Sigma A^{1/2}G'= ([I_{k}, 0] \vdots[0, 0]
 )'$. Define $V=GA^{1/2}X$. Then,
$\E[h(X'AX)WA^{1/2\dag}A^{1/2}X]=\E[h(V_1'V_1)WA^{1/2\dag}
G'[I_k,0]'V_1]$ with $V_{1}=[I_k,0]GA^{1/2}V$, and then, the rest of
the proof follows from Theorem~1 in Judge and Bock \cite
{judg1978} along with
some algebraic computations.
\end{pf}

%
\begin{rem}
For the special case where $\Sigma$ is the $p$-dimensional
 matrix $I_{p}$, Theorem~\ref{thma1} gives Theorem~1 in
Judge and Bock \cite{judg1978} with $A=W^{*}=I_{p}$. This
shows that the
provided theorem generalizes the quoted classical result.
\end{rem}

By using Theorem~\ref{thma1}, we establish the following corollary.

%
\begin{cor}\label{prna1}
Set $\mu_2=-\mu_1$ and let
$\epsilon_5$ be as defined in Lemma~\ref{norm1}.
Let $h$ be a Borel measurable and real-valued integrable function,
let $W=A^{1/2}W^*A^{1/2}$, $W^*$ is a nonnegative definite matrix.
Then, we have $\E[h(\epsilon_5'A\epsilon_5)W\epsilon_5]=\E
[h(\chi_{k+2}^2(\mu_2'A\mu_2))]W\mu_2$.
\end{cor}

%
\begin{thm}\label{thmb1}
Let $D_1=\trace(W\Sigma)$, $D_2=\mu'W\mu$ and assume the
conditions of Theorem~\textup{\ref{thma1}} hold. Then, $\E[h(X'AX)X'WX]=\E
[h(\chi_{k+2}^2(\mu'A\mu))]D_1
+\E[h(\chi_{k+4}^2(\mu'A\mu))]D_2$.
\end{thm}

\begin{pf}
By using the same transformation methods as in the proof of
Theorem~\ref{thma1}, we have
\begin{eqnarray*}
\E\bigl[h\bigl(X'AX\bigr)X'WX\bigr]=\E\bigl[h
\bigl(V_1'V_1\bigr)V_1'[I_k,0]GA^{1/2\dag}
WA^{1/2\dag}G'[I_k,0]'V_1
\bigr].
\end{eqnarray*}
Therefore, the proof is completed by combining Theorem~2 in Judge and Bock \cite{judg1978} along with some algebraic computations.
\end{pf}

%
\begin{rem}
Note that Theorem~\ref{thmb1} generalizes Theorem~2 in Judge and Bock \cite{judg1978}. Indeed, if $\Sigma=I_{p}$, the quoted
result is
obtained by taking $A=I_{p}$.
\end{rem}

By using Theorem~\ref{thmb1}, we establish the following corollary.

%
\begin{cor}\label{prnb1}
Let $D_1=\trace(W\Lambda_{11})$,
$D_2=\mu_2'W\mu_2$ and suppose that the conditions of Corollary~\textup{\ref
{prna1}} hold. Then,
$\E[h(\epsilon_5'A\epsilon_5)\epsilon_5'W\epsilon_5]
=\E[h(\chi_{k+2}^2(\mu_2'A\mu_2))]D_1 +\E[h(\chi_{k+4}^2(\mu
_2'A\mu_2))]D_2$.
\end{cor}

\begin{pf}This corollary directly follows from Theorem~\ref{thmb1}.
\end{pf}

%
\begin{thm}\label{lema2}
Let
\[
\lleft(
\begin{array}{c}
X\\
Y
\end{array}
 \rright)\sim\mathcal{N}_{2p}\lleft(
 \left(
\begin{array}{c}
\mu_{X}\\
\mu_{Y}
\end{array}
\right ),
 \left(
\begin{array}{c@{\quad}c}
\Sigma_{11}&\Sigma_{12}\\
\Sigma_{21}&\Sigma_{22}
\end{array}
 \right)
 \rright),
 \]
 where the rank of $\Sigma_{11}$ is $k\leqslant p$, with
$\mu_Y=-\mu_X$, $A\Sigma_{11}A=A$;
$\Sigma_{11}A\Sigma_{11}=\Sigma_{11}$; $\Sigma_{11}A\mu_X=\mu_X$.
Further, we assume that $W=A^{1/2}W^*A^{1/2}$, where $W^*$ is a
nonnegative definite matrix. Then,
\begin{eqnarray*}
&&\E\bigl[h\bigl(X'AX\bigr)Y'WX\bigr]\\
&&\quad =-\E\bigl[h\bigl(
\chi_{k+2}^2\bigl(\mu_X'A
\mu_X\bigr)\bigr)\bigr]\mu_X'W
\mu_X -\E\bigl[h\bigl(\chi_{k+2}^2\bigl(
\mu_X'A\mu_X\bigr)\bigr)\bigr]
\mu_X'A\Sigma _{12}W\mu_X
\\
&&\qquad {}+\E\bigl[h\bigl(\chi_{k+2}^2\bigl(\mu_X'A
\mu_X\bigr)\bigr)\bigr]\trace(\Lambda _{12}W
\Lambda_{11}A) +\E\bigl[h\bigl(\chi_{k+4}^2\bigl(
\mu_X'A\mu_X\bigr)\bigr)\bigr]
\mu_X'A\Lambda_{12}W\mu_X.
\end{eqnarray*}
\end{thm}

\begin{pf}
Using the similar transformation methods as in proof of
Theorem~\ref{thma1}, we have
\begin{eqnarray*}
\E\bigl[h\bigl(X'AX\bigr)Y'WX\bigr] =\E\bigl[h
\bigl(V_1'V_1\bigr)\E[Y|V_1]'WA^{1/2\dag}G'[I_k,0]'V_1
\bigr],
\end{eqnarray*}
where
$\E[Y|V_1]=-\mu_X+\Sigma_{21}A^{1/2}G'[I_k,0]'(V_1-\mu_v)$.
Further, from Theorem~\ref{thma1},
\begin{eqnarray*}
\E\bigl[h\bigl(V_1'V_1\bigr)
\mu_2'WA^{1/2\dag}G'[I_k,0]'V_1
\bigr]  =\E\bigl[h\bigl(\chi_{k+2}^2\bigl(
\mu_X'A\mu_X\bigr)\bigr)\bigr]
\mu_X'W\mu_X
\end{eqnarray*}
and
\begin{eqnarray*}
&&\E\bigl[h\bigl(V_1'V_1\bigr)
\mu_v'[I_k,0]GA^{1/2}
\Sigma_{12}WA^{1/2\dag}G'[I_k,0]'V_1
\bigr]\\
&&\quad  =\E\bigl[h\bigl(\chi_{k+2}^2\bigl(
\mu_X'A\mu_X\bigr)\bigr)\bigr]
\mu_X'A\Sigma_{12}WA^{1/2\dag
}
\mu_X,
\end{eqnarray*}
and the proof is completed by some algebraic computations.
\end{pf}

By using this theorem, we establish the following corollary.

%
\begin{cor}
With $\epsilon_5$ and $\epsilon_4$ defined in Lemma~\ref{norm1}, and
let $\mu_2=-\mu_1$.
Then, we have
\begin{eqnarray*}
&&\mathrm{E}\bigl[h\bigl(\epsilon_5'A
\epsilon_5\bigr)\epsilon_4'W
\epsilon_5\bigr]
\\
&&\quad =-\mathrm{E}\bigl[h\bigl(\chi_{k+2}^2\bigl(
\mu_2'A\mu_2\bigr)\bigr)\bigr]
\mu_2'W\Lambda _{11}A\mu_2 -
\mathrm{E}\bigl[h\bigl(\chi_{k+2}^2\bigl(
\mu_2'A\mu_2\bigr)\bigr)\bigr]
\mu_2'A\Lambda _{12}W\Lambda_{11}A
\mu_2
\\
&&\qquad {}+\mathrm{E}\bigl[h\bigl(\chi_{k+2}^2\bigl(
\mu_2'A\mu_2\bigr)\bigr)\bigr]\trace(
\Lambda _{12}W\Lambda_{11}A)\\
&&\qquad {} +\mathrm{E}\bigl[h\bigl(
\chi_{k+4}^2\bigl(\mu_2'A
\mu_2\bigr)\bigr)\bigr]\mu_2'A\Lambda
_{12}W\Lambda_{11}A\mu_2.
\end{eqnarray*}
\end{cor}

\begin{pf*}{Proof of Corollary~\ref{adrcompar}}
By some algebraic computations, we have,
\begin{eqnarray*}
&&\ADR\bigl(\hat{\delta}^s,\delta^0,W\bigr)-\ADR\bigl(
\hat{\delta},\delta ^0,W\bigr)
\\
&&\quad =-(k-2)^2\trace\bigl(W(\Lambda_{11}+2
\Lambda_{12})\bigr)\E\bigl[\chi _{k+2}^{-4}(\Delta)
\bigr]\\
&&\qquad {}-(k-2) \bigl(4\Delta C_1-(k+2)C_2\bigr)\E\bigl[
\chi_{k+4}^{-4}(\Delta)\bigr],
\end{eqnarray*}
where $C_1=\trace(W(\Lambda_{11}+\Lambda_{12}))$,
$C_2=\mu_{1}^{\prime} A(\Lambda_{11} +4\Lambda_{12}/(k+2))W\mu_1$, and
$C_3=\linebreak[4] \trace(W\Lambda_{11})$. Then, since $k\geqslant2$,
$\ADR(\hat{\delta}^s,\delta^0,W)\leqslant
\ADR(\hat{\delta},\delta^0,W)$ provided that\linebreak[4]
$\trace(W(\Lambda_{11}+2\Lambda_{12}))\geqslant0$ and
$4\Delta C_1-(k+2)C_2\geqslant0$.
Note that if $C_2=0$, $4\Delta C_1-(k+2)C_2\geqslant0$ holds for
any $\Delta\geqslant0$, and if $C_2>0$, $4\Delta
C_1-(k+2)C_2\geqslant0$ holds for $\Delta C_1\geqslant(k+2)C_2/4$,
which is equivalent to $C_1\geqslant(k+2)C_2/(4\Delta)$.

Since $C_2=\mu'_1
A(\Lambda_{11}+4\Lambda_{12}/(k+2))W\Lambda_{11}A\mu_1$, and by
Courant's theorem, we have
\[
\Ch_{\mathrm{min}}\bigl(\Pi^*\bigr)\leqslant\frac{\mu_1^{\prime}
A(\Lambda_{11}+4\Lambda_{12}/(k+2)) W\Lambda_{11}A\mu_1}{\mu
_1^{\prime}
A\mu_1} \leqslant
\Ch_{\max}\bigl(\Pi^*\bigr),
\]
where $\Pi^*=(\Pi_{0}+\Pi'_{0})/2$,
$\Pi_{0}=A^{1/2}(\Lambda_{11}+4\Lambda_{12}/(k+2))W\Lambda_{11}A^{1/2}$
and $\Ch_{\mathrm{min}}(\Pi^*)$, $\Ch_{\max}(\Pi^*)$ are
denoted as the smallest and largest eigenvalue of $\Pi^*$,
respectively. Then, $4\Delta C_1-(k+2)C_2\geqslant0$ holds if
$C_1\geqslant(k+2)\Ch_{\max}(\Pi^*)/4$. In addition,
since
$\trace(W(\Lambda_{11}+2\Lambda_{12}))\geqslant0$ is
equivalent to $C_1\geqslant-\trace(W\Lambda_{12})$, it
follows that
\[
\ADR\bigl(\hat{\delta}^s, \delta^0,W\bigr)\leqslant \ADR
\bigl(\hat{\delta},\delta^0,W\bigr)
\]
if
$\trace(W(\Lambda_{11}+\Lambda_{12}))\geqslant
\max(-\trace(W\Lambda_{12}),(k+2)\Ch_{\max}(\Pi^*)/4)$.
Further, by similar algebraic computations, we prove that
$\ADR(\hat{\delta}^{s+},\delta^0,W)
\leqslant\ADR(\hat{\delta}^s,\delta^0,W)$, this completes
the proof.
\end{pf*}
\end{appendix}

\section*{Acknowledgements}
The authors would like to acknowledge the financial support received
from Natural Sciences and Engineering Research Council of Canada.
Further, the authors would like to thank anonymous referees for
useful comments and suggestions.

%

\printhistory

\begin{thebibliography}{20}


\bibitem{bp2003}
\begin{barticle}[auto:STB|2014/06/18|12:29:53]
\bauthor{\bsnm{Bai},~\bfnm{J.}\binits{J.}} \AND
\bauthor{\bsnm{Perron},~\bfnm{P.}\binits{P.}}
(\byear{2003}).
\btitle{Computation and analysis of multiple structural change models}.
\bjournal{J. Appl. Econometr.}
\bvolume{18}
\bpages{1--22}.
\end{barticle}
\bptok{imsref}%
\endbibitem

\bibitem{baranchick1964}
\begin{bmisc}[auto:STB|2014/06/18|12:29:53]
\bauthor{\bsnm{Baranchick},~\bfnm{A.}\binits{A.}}
\bhowpublished{(1964). Multiple regression and estimation of the mean of a multivariate normal distribution.
Technical Report No. 51, Dept.  Statistics, Stanford Univ.}
\end{bmisc}
\bptok{imsref}%
\endbibitem

\bibitem{bm1998}
\begin{barticle}[auto:STB|2014/06/18|12:29:53]
\bauthor{\bsnm{Braun},~\bfnm{J.~V.}\binits{J.V.}} \AND
\bauthor{\bsnm{Muller},~\bfnm{H.~G.}\binits{H.G.}}
(\byear{1998}).
\btitle{Statistical methods for DNA sequence segmentation}.
\bjournal{Statist. Sci.}
\bvolume{13}
\bpages{142--162}.
\end{barticle}
\bptok{imsref}%
\endbibitem

\bibitem{bt1987}
\begin{bbook}[mr]
\bauthor{\bsnm{Broemeling},~\bfnm{Lyle~D.}\binits{L.D.}} \AND
\bauthor{\bsnm{Tsurumi},~\bfnm{Hiroki}\binits{H.}}
(\byear{1987}).
\btitle{Econometrics and Structural Change}.
\bseries{Statistics: Textbooks and Monographs}
\bvolume{74}.
\blocation{New York}:
\bpublisher{Dekker, Inc.}
\bid{mr={0922263}}
\end{bbook}
\bptok{imsref}%
\endbibitem

\bibitem{fc1990a}
\begin{barticle}[mr]
\bauthor{\bsnm{Fu},~\bfnm{Yun-Xin}\binits{Y.-X.}} \AND
\bauthor{\bsnm{Curnow},~\bfnm{R.~N.}\binits{R.N.}}
(\byear{1990}).
\btitle{Locating a changed segment in a sequence of {B}ernoulli variables}.
\bjournal{Biometrika}
\bvolume{77}
\bpages{295--304}.
\bid{doi={10.1093/biomet/77.2.295}, issn={0006-3444}, mr={1064801}}
\end{barticle}
\bptok{imsref}%
\endbibitem

\bibitem{fc1990b}
\begin{barticle}[mr]
\bauthor{\bsnm{Fu},~\bfnm{Yun-Xin}\binits{Y.-X.}} \AND
\bauthor{\bsnm{Curnow},~\bfnm{R.~N.}\binits{R.N.}}
(\byear{1990}).
\btitle{Maximum likelihood estimation of multiple change points}.
\bjournal{Biometrika}
\bvolume{77}
\bpages{563--573}.
\bid{doi={10.1093/biomet/77.3.563}, issn={0006-3444}, mr={1087847}}
\end{barticle}
\bptok{imsref}%
\endbibitem

\bibitem{ahm2009}
\begin{barticle}[mr]
\bauthor{\bsnm{Hossain},~\bfnm{S.}\binits{S.}},
\bauthor{\bsnm{Doksum},~\bfnm{Kjell~A.}\binits{K.A.}} \AND
\bauthor{\bsnm{Ahmed},~\bfnm{S.E.}\binits{S.E.}}
(\byear{2009}).
\btitle{Positive shrinkage, improved pretest and absolute penalty estimators in partially linear models}.
\bjournal{Linear Algebra Appl.}
\bvolume{430}
\bpages{2749--2761}.
\bid{doi={10.1016/j.laa.2008.12.015}, issn={0024-3795}, mr={2509855}}
\end{barticle}
\bptok{imsref}%
\endbibitem

\bibitem{js1961}
\begin{bincollection}[mr]
\bauthor{\bsnm{James},~\bfnm{W.}\binits{W.}} \AND
\bauthor{\bsnm{Stein},~\bfnm{Charles}\binits{C.}}
(\byear{1961}).
\btitle{Estimation with quadratic loss}.
In \bbooktitle{Proc. 4th {B}erkeley {S}ympos. {M}ath. {S}tatist. and {P}rob.}
\bvolume{I}
\bpages{361--379}.
\blocation{Berkeley, CA}:
\bpublisher{Univ. California Press}.
\bid{mr={0133191}}
\end{bincollection}
\bptok{imsref}%
\endbibitem

\bibitem{judg1978}
\begin{bbook}[mr]
\bauthor{\bsnm{Judge},~\bfnm{George~G.}\binits{G.G.}} \AND
\bauthor{\bsnm{Bock},~\bfnm{M.~E.}\binits{M.E.}}
(\byear{1978}).
\btitle{The Statistical Implications of Pre-Test and {S}tein-Rule Estimators in Econometrics}.
\blocation{Amsterdam}:
\bpublisher{North-Holland}.
\bid{mr={0483199}}
\end{bbook}
\bptok{imsref}%
\endbibitem

\bibitem{judg2004}
\begin{barticle}[mr]
\bauthor{\bsnm{Judge},~\bfnm{George~G.}\binits{G.G.}} \AND
\bauthor{\bsnm{Mittelhammer},~\bfnm{Ron~C.}\binits{R.C.}}
(\byear{2004}).
\btitle{A semiparametric basis for combining estimation problems under quadratic loss}.
\bjournal{J. Amer. Statist. Assoc.}
\bvolume{99}
\bpages{479--487}.
\bid{doi={10.1198/016214504000000430}, issn={0162-1459}, mr={2062833}}
\end{barticle}
\bptok{imsref}%
\endbibitem

\bibitem{l1986}
\begin{barticle}[auto:STB|2014/06/18|12:29:53]
\bauthor{\bsnm{Lombard},~\bfnm{F.}\binits{F.}}
(\byear{1986}).
\btitle{The change-point problem for angular data: A nonparametric approach}.
\bjournal{Technometrics}
\bvolume{28}
\bpages{391--397}.
\end{barticle}
\bptok{imsref}%
\endbibitem

\bibitem{mcleish1977}
\begin{barticle}[mr]
\bauthor{\bsnm{McLeish},~\bfnm{D.~L.}\binits{D.L.}}
(\byear{1977}).
\btitle{On the invariance principle for nonstationary mixingales}.
\bjournal{Ann. Probab.}
\bvolume{5}
\bpages{616--621}.
\bid{mr={0445583}}
\end{barticle}
\bptok{imsref}%
\endbibitem

\bibitem{n2011}
\begin{barticle}[mr]
\bauthor{\bsnm{Nkurunziza},~\bfnm{S{\'e}v{\'e}rien}\binits{S.}}
(\byear{2011}).
\btitle{Shrinkage strategy in stratified random sample subject to measurement error}.
\bjournal{Statist. Probab. Lett.}
\bvolume{81}
\bpages{317--325}.
\bid{doi={10.1016/j.spl.2010.10.020}, issn={0167-7152}, mr={2764300}}
\end{barticle}
\bptok{imsref}%
\endbibitem

\bibitem{n2012a}
\begin{barticle}[mr]
\bauthor{\bsnm{Nkurunziza},~\bfnm{S{\'e}v{\'e}rien}\binits{S.}}
(\byear{2012}).
\btitle{The risk of pretest and shrinkage estimators}.
\bjournal{Statistics}
\bvolume{46}
\bpages{305--312}.
\bid{doi={10.1080/02331888.2010.508561}, issn={0233-1888}, mr={2929155}}
\end{barticle}
\bptok{imsref}%
\endbibitem

\bibitem{nku}
\begin{barticle}[mr]
\bauthor{\bsnm{Nkurunziza},~\bfnm{S{\'e}v{\'e}rien}\binits{S.}} \AND
\bauthor{\bsnm{Ahmed},~\bfnm{S.~Ejaz}\binits{S.E.}}
(\byear{2010}).
\btitle{Shrinkage drift parameter estimation for multi-factor {O}rnstein--{U}hlenbeck processes}.
\bjournal{Appl. Stoch. Models Bus. Ind.}
\bvolume{26}
\bpages{103--124}.
\bid{doi={10.1002/asmb.775}, issn={1524-1904}, mr={2722886}}
\end{barticle}
\bptok{imsref}%
\endbibitem

\bibitem{pq2006}
\begin{barticle}[mr]
\bauthor{\bsnm{Perron},~\bfnm{Pierre}\binits{P.}} \AND
\bauthor{\bsnm{Qu},~\bfnm{Zhongjun}\binits{Z.}}
(\byear{2006}).
\btitle{Estimating restricted structural change models}.
\bjournal{J. Econometrics}
\bvolume{134}
\bpages{373--399}.
\bid{doi={10.1016/j.jeconom.2005.06.030}, issn={0304-4076}, mr={2328414}}
\end{barticle}
\bptok{imsref}%
\endbibitem

\bibitem{by2009}
\begin{barticle}[mr]
\bauthor{\bsnm{Perron},~\bfnm{Pierre}\binits{P.}} \AND
\bauthor{\bsnm{Yabu},~\bfnm{Tomoyoshi}\binits{T.}}
(\byear{2009}).
\btitle{Testing for shifts in trend with an integrated or stationary noise component}.
\bjournal{J. Bus. Econom. Statist.}
\bvolume{27}
\bpages{369--396}.
\bid{doi={10.1198/jbes.2009.07268}, issn={0735-0015}, mr={2554242}}
\end{barticle}
\bptok{imsref}%
\endbibitem

\bibitem{sal}
\begin{bbook}[mr]
\bauthor{\bsnm{Saleh},~\bfnm{A.~K.~Md.~Ehsanes}\binits{A.K.Md.E.}}
(\byear{2006}).
\btitle{Theory of Preliminary Test and {S}tein-Type Estimation with Applications}.
\bseries{Wiley Series in Probability and Statistics}.
\blocation{Hoboken, NJ}:
\bpublisher{Wiley}.
\bid{doi={10.1002/0471773751}, mr={2218139}}
\end{bbook}
\bptok{imsref}%
\endbibitem

\bibitem{tan}
\begin{bmisc}[auto:STB|2014/06/18|12:29:53]
\bauthor{\bsnm{Tan},~\bfnm{Z.}\binits{Z.}}
\bhowpublished{(2014). Improved minimax estimation of a multivariate normal mean under heteroscedasticity. \textit{Bernoulli}. To appear}.
\end{bmisc}
\bptok{imsref}%
\endbibitem

\bibitem{zei}
\begin{barticle}[mr]
\bauthor{\bsnm{Zeileis},~\bfnm{Achim}\binits{A.}},
\bauthor{\bsnm{Kleiber},~\bfnm{Christian}\binits{C.}},
\bauthor{\bsnm{Kr{\"a}mer},~\bfnm{Walter}\binits{W.}} \AND
\bauthor{\bsnm{Hornik},~\bfnm{Kurt}\binits{K.}}
(\byear{2003}).
\btitle{Testing and dating of structural changes in practice}.
\bjournal{Comput. Statist. Data Anal.}
\bvolume{44}
\bpages{109--123}.
\bid{doi={10.1016/S0167-9473(03)00030-6}, issn={0167-9473}, mr={2019790}}
\end{barticle}
\bptok{imsref}%
\endbibitem
\end{thebibliography}
\end{document}